\documentclass[sigconf,authorversion]{acmart}

\usepackage{balance} % For balanced columns on the last page
\usepackage{booktabs} % For formal tables
\usepackage{array} % for better arrays (eg matrices) in maths
\usepackage{dsfont}
\usepackage{algorithm}
\usepackage{algorithmicx}
\usepackage{algpseudocode}
\usepackage{url}
\usepackage{pythonhighlight}
\usepackage{xspace}

\newcommand{\SUMO}{{\textrm{Sofomore}}\xspace}
\newcommand{\COMOCMA}{\textrm{COMO-CMA-ES}\xspace} % might change

\DeclareMathOperator*{\argmax}{arg\,max}

\newcommand{\HVC}{\ensuremath{\text{HVC}}\xspace}
\newcommand{\HVI}{\ensuremath{\text{HVI}}\xspace}
\newcommand{\UHVI}{\ensuremath{\text{UHVI}}\xspace}

\newcommand{\f}{{\mathbf f}}
\newcommand{\kernel}{kernel\xspace}
\newcommand{\kernels}{kernels\xspace}
\newcommand{\convergencegap}{convergence gap\xspace}
\newcommand{\archivegap}{archive gap\xspace}

\makeatletter

\newcommand{\Rmnum}[1]{\expandafter\@slowromancap\romannumeral #1@} % roman numerals : I, II, ...
\makeatother

\newcommand{\markupdraft}[2]{% {#1: {color|display} command}{#2: desired color or text}
    \ifthenelse{\equal{#1}{display}}{#2}{}%                 % display only in draft version
    \ifthenelse{\equal{#1}{color}}{\color{#2}}{}%           % colored only in draft (for \new command)
}
\newcommand{\nnotecolored}[3][]{\markupdraft{display}{{\color{#2}#3$_{-\!\mathrm{#1}}$}}}
\newcommand{\newcolored}[3][]{{\markupdraft{color}{#2}#3}%  % kept in the final print
\ifthenelse{\equal{#1}{}}{}{\markupdraft{display}{{\color{yellow!70!black}[#1]}}}}

\newcommand{\del}[2][]{{\markupdraft{display}{{\color{orange}[removed: "#2"[#1]]}}}} % (to be) removed
\newcommand{\new}[2][]{\newcolored[#1]{blue!90!black}{#2}}%  % kept in the final print
\newcommand{\nnew}[2][]{\newcolored[#1]{red}{#2}}%  % kept in the final print
\renewcommand{\del}[2][]{}  % do not show "deleted" text anymore

\newcommand{\niko}[1]{\protect\nnotecolored[Niko]{green!50!black}{~#1}}

\newcommand{\ba}{\begin{eqnarray}}
\newcommand{\ea}{\end{eqnarray}}
\newcommand{\baStar}{\begin{eqnarray*}}
\newcommand{\eaStar}{\end{eqnarray*}}

     \newcommand{\foncfast}[4]{#1 \ni #2 \longmapsto #3 \in #4}
  \newcommand{\fonccourte}[2]{#1 \longmapsto #2}

\newcommand{\dsp}{\displaystyle}
\newcommand{\acco}[1]{\left\{#1\right\}} 
\newcommand{\pare}[1]{\left(#1\right)}
\newcommand{\croc}[1]{\left[#1\right]}

\DeclareMathOperator{\mydef}{def}

\newcommand{\egaldef}{\stackrel{{\scriptscriptstyle{\mydef}}}{=}}

\def\R{{\mathbb R}} \def\diff{{\mathrm d}}
\def\N{{\mathbb N}}  \def\F{\mathcal{F}} \def\E{\mathcal{E}} \def\I{\mathcal{I}}\def\T{\Theta} 
\def\emp{ {\mathbf E\mathbf P\mathbf F} }

\def\r{\mathbf{r}}
\def\dimnumber{n}
\def\kernumber{p}
\def\objnumber{k}

\definecolor{grey}{rgb}{0.95,0.95,0.95}
\renewcommand{\markupdraft}[2]{}  % final mode text command

\begin{document}

\copyrightyear{2019} 
\acmYear{2019} 
\setcopyright{acmlicensed}
\acmConference[GECCO '19]{Genetic and Evolutionary Computation Conference}{July 13--17, 2019}{Prague, Czech Republic}
\acmBooktitle{Genetic and Evolutionary Computation Conference (GECCO '19), July 13--17, 2019, Prague, Czech Republic}
\acmPrice{15.00}
\acmDOI{10.1145/3321707.3321852}
\acmISBN{978-1-4503-6111-8/19/07}

\title[\COMOCMA and the \SUMO framework]{
%COMO-CMA-ES: A Linearly Convergent Multiobjective Algorithm
%Single-objective Optimization of Multiobjective Problems Via Indicator-Based Subspace Optimization: The Sofomore Framework and the COMO-CMA-ES
%A framework for defining new multiobjective algorithms from single-objective ones  Instantiated on CMA-ES: the COMO-CMA-ES
%COMO-CMA-ES: A fast COmma Multiobjective Covariance Matrix Adaptation}
%\title{Indicator-Based Subspace Optimization of Multiobjective Problems: \COMOCMA and the \SUMO framework}
%Penalized Hypervolume Improvement for Multiobjective Problems: \COMOCMA and the \SUMO framework
Uncrowded Hypervolume Improvement: \COMOCMA and the \SUMO framework}

\newcommand{\ouraffiliation}{
\affiliation{%
  \institution{Inria, CMAP, Ecole Polytechnique, IP Paris}
%  \streetaddress{}
% \city{Paris Area} 
%  \state{France} 
%  \postcode{}
}}
\newcommand{\shortaffiliation}{
\affiliation{%
  \institution{Inria}%, CMAP, Ecole Polytechnique, IP Paris}
%  \streetaddress{}
% \city{Paris Area} 
%  \state{France} 
%  \postcode{}
}}

%% The submitted version for review should be ANONYMOUS
\author{Cheikh Tour\'{e}}%$^{\star}$}
%\authornote{Cheikh insisted to be corresponding author}
\ouraffiliation
\email{cheikh.toure@polytechnique.edu}
\author{Nikolaus Hansen}
\ouraffiliation
\email{firstname.lastname@inria.fr}
\author{Anne Auger}
\ouraffiliation
\email{firstname.lastname@inria.fr}
\author{Dimo Brockhoff}
\ouraffiliation
\email{firstname.lastname@inria.fr}

%\thanks{$^\star$Corresponding author.}% Please use \url{cheikh.toure@polytechnique.edu}.}

% The default list of authors is too long for heaaders.
\renewcommand{\shortauthors}{C.\ Tour\'{e} et al.}

\begin{abstract}
We present a framework to build a multiobjective algorithm from single-objective ones.\del{Denoting $n$ the search space dimension, }
This framework \nnew{addresses}\del{approaches} the $p \times n$-dimen\-sional problem of finding $p$ solutions \new{in an $n$-dimensional search space,} maximizing an indicator by dynamic subspace optimization. Each single-objective algorithm optimizes the indicator function given $p-1$ fixed solutions. Crucially, dominated solutions minimize their distance to the empirical Pareto front defined by these $p-1$ solutions. We instantiate the framework with CMA-ES as single-objective optimizer.
The new algorithm, \COMOCMA, is empirically shown to converge linearly on bi-objective convex-quadratic problems and is compared to MO-CMA-ES, NSGA-II and SMS-EMOA.
\end{abstract}



\maketitle

%\input{outline-gecco-paper.tex}

%%%%%%%%%%%%%%%%%%%%%%%%%%%%%%%%%%%%%%%%%%%%%%%%%%%%%%%%%%%%%%%%%%%%%
\section{Introduction}
%%%%%%%%%%%%%%%%%%%%%%%%%%%%%%%%%%%%%%%%%%%%%%%%%%%%%%%%%%%%%%%%%%%%%

Multiobjective optimization problems must be solved frequently in practice. In contrast to the optimization of a single objective, solving a multiobjective problem involves to handle trade-offs or incomparabilities between the objective functions such that the aim is to approximate the Pareto set---the set of all Pareto-optimal or non-dominated solutions. One might be interested to obtain an approximation of unbounded size (the more points the better) or just to have $\kernumber$ points approximating the Pareto set.
Evolutionary Multiobjective Optimization (EMO) algorithms aim at such an approximation in a single algorithm run whereas more classical approaches, e.g.\ optimizing a weighted sum of the objectives with changing weights, operate in multiple runs.

The first introduced EMO algorithms simply changed the selection of an existing single-objective evolutionary algorithm keeping the exact same search operators. The population at a given iteration was then providing an approximation of the Pareto set. This idea led to the practically highly successful NSGA-II algorithm \cite{dapm2002a} that employs a two-step fitness assignment: after a first non-dominated ranking \cite{gold1989a}, solutions with equal non-domination rank are further distinguished by their crowding distance---based on the distance of each solution to its neighbors in objective space. However, it has been pointed out that NSGA-II does not converge to the Pareto set in a mathematical sense due to so-called deteriorative cycles: if all population members of the algorithm are non-dominated at some point in time, it is only the crowding distance that is optimized, without indicating any search direction towards the Pareto set to the algorithm. As a result, solutions which had been non-dominated solutions at some point in time can be replaced by previously dominated ones during the optimization, ending up in a cyclic but not in convergent behavior \cite{bfn2010a}.

To improve the convergence properties of EMO algorithms, different approaches have been introduced later, most notably the indicator-based algorithms and especially algorithms based on the hypervolume indicator. They replace the crowding distance of NSGA-II with the (hypervolume) indicator contribution, see e.g.\ \cite{ihr2007a,bne2007a}. Using the hypervolume indicator has the advantage that it is the only known strictly monotone quality indicator \cite{ktz2006a} (see also next section) and thus, its optimization will result in solution sets that are subsets of the Pareto set.

The optimization goal of indicator-based algorithms such as SMS-EMOA \cite{bne2007a} or MO-CMA-ES \cite{ihr2007a} is to find the best set of $\kernumber$ solutions with respect to a given quality indicator (the set with the largest quality indicator value among all sets of size $\kernumber$). This optimal set of $\kernumber$ solutions is known as the optimal $\kernumber$-distribution~\cite{auger2009theory}. In principle, the search for the optimal $\kernumber$-distribution can be formalized as a $\kernumber\cdot\dimnumber$-dimensional optimization problem where $\kernumber$ is the number of solutions and $\dimnumber$ is the dimension of the search space. 

% \niko{From here on the introduction is an narrated table of contents, but this is a probably good preparation of the reader for the more technical parts.}
As we will discuss later, it turns out that this optimization problem is not only of too high dimension in practice but also flat in large regions of the search space if the hypervolume indicator is the underlying quality indicator. The combination of non-dominated ranking and hypervolume contribution as in SMS-EMOA or MO-CMA-ES corrects for this flatness, but also introduces search directions that are pointing towards already existing non-dominated solutions and not towards not-yet-covered regions of the Pareto set.
In this paper, we show that we can correct the flat region of the hypervolume indicator by introducing a search bias towards yet-uncovered regions of the Pareto set by adding the distance to the empirical non-domination front, which leads to the new notion of Uncrowded Hypervolume Improvement. Then, we define a (dynamic) fitness function that can be optimized by single-objective algorithms. From there, going back to this original idea of EMO algorithms to use single-objective optimizers to build an EMO,
we define the \emph{Single-objective Optimization FOr Optimizing Multiobjective Optimization pRoblEms} framework (\SUMO) to build in an elegant manner, a multiobjective algorithm from a set of $\kernumber$ single-objective optimizers. Each single-objective algorithm optimizes (iteratively or in parallel) a dynamic fitness that depends on the output of the other $\kernumber-1$ optimizers. 
%This differs from the idea of MOEA/D in which $\kernumber$ single-objective optimizers optimize in parallel their own fixed fitness (``scalarizing'') function which are coupled by exchanging solutions.

We instantiate the \SUMO framework with the state-of-the-art single-objective algorithm CMA-ES. We show experimentally that the ensuing \COMOCMA (Comma-Selection Multiobjective CMA-ES) exhibits linear convergence towards the optimal $\kernumber$-distribution on a wide variety of bi-objective convex quadratic functions. In contrast, default implementations of the SMS-EMOA where the reference point is\del{ not} fixed and NSGA-II do not exhibit this linear convergence. The comparison between \COMOCMA and a previous MATLAB implementation of the elitist MO-CMA-ES also shows the same or an improved convergence \emph{speed} in \COMOCMA except for the double sphere function.

The paper is structured as follows. In the next section, we start with preliminaries related to multiobjective optimization and quality indicators. Section~\ref{sec:sofomore} discusses the fitness landscape of indicator- and especially hypervolume-based quality measures and eventually introduces our \SUMO framework. Section~\ref{sec:comocma} gives details about the new \COMOCMA algorithm as an instantiation of \SUMO with CMA-ES. Section~\ref{sec:experiments} experimentally validates the new algorithm and compares it with three existing algorithms from the literature and Section~\ref{sec:conclusions} discusses the results and concludes the paper.

%%%%%%%%%%%%%%%%%%%%%%%%%%%%%%%%%%%%%%%%%%%%%%%%%%%%%%%%%%%%%%%%%%%%%
\section{Preliminaries}\label{sec:prelim}
%%%%%%%%%%%%%%%%%%%%%%%%%%%%%%%%%%%%%%%%%%%%%%%%%%%%%%%%%%%%%%%%%%%%%
In the following, we assume without loss of generality the minimization of a vector-valued function
%\begin{equation}
$\f: x \in \R^\dimnumber \mapsto \f(x)=(\f_1(x),\ldots,\f_\objnumber(x)) \in \R^\objnumber$
%\end{equation}
that maps a search point from the search space $\R^{\dimnumber}$ of dimension $\dimnumber$
to the objective space $\f(\R^{\dimnumber}) \subseteq \R^\objnumber$. 
This minimization of $\f$ is generally formalized in terms of the weak Pareto dominance
relation for which we write that a search point $x\in\R^{\dimnumber}$
weakly Pareto-dominates another search point $y\in\R^{\dimnumber}$ (written in short as $x \preceq y$ or with an abuse of notation as $\f(x) \preceq \f(y)$) if and only if $\f_{i}(x) \leq \f_{i}(y)$ for all $i \in \{1, \ldots, \objnumber\}$. Note also that we can naturally extend the (weak) Pareto dominance relation to \del{(finite) }subsets $A,B\subset \R^{\dimnumber}$ as $A \preceq B$ if and only if for all $ b\in B$, there exists $a\in A$ such that $a \preceq b$.
If the relation $\leq$ is strict for at least one objective function, we say that $x$ Pareto-dominates $y$ (and write $x \prec y$).
The set of non-dominated search points constitutes the so-called Pareto set, its image under $\f$ is called the Pareto front.
In the remainder, we will also use the term \emph{empirical non-dominated front} or  \emph{empirical Pareto front} ($\emp_{S,\r}$) for objective vectors that are on the boundary of the (objective space) region dominating a reference point $\r\in\R^{\objnumber}$, and not dominated by any element of $\f(S)$ with $S\subset \R^{\dimnumber}$: 
\begin{equation}
\label{eq:empiricalnondominationfront}
  \emp_{S,\r} = \partial U_{S,\r}, \text{ with }  U_{S,\r} = \acco{ z \prec \r; \forall s \in S, \f(s) \nprec z } 
\end{equation}
where $\partial U_{S,\r}$ is the boundary of the non-dominated region $U_{S,\r}$. Note that $\emp_{S,\r}\cap\f(\R^\dimnumber)$ is the Pareto front when $S$ contains the Pareto set.
\paragraph{Indicator-Based Set Optimization Problems}
Pareto sets and Pareto fronts are, under mild assumptions, $\objnumber-1$ dimensional manifolds. In practice, we are often interested in a finite size approximation of these sets with, let us say, $\kernumber \,(\geq 1)$ many search points. 
To assess the quality of a Pareto set approximation $S\subseteq \R^{\dimnumber}$, a quality indicator $\I: 2^{\R^{\dimnumber}} \rightarrow \R$ assigns a real valued quality $I(S)$ to $S$. Formally speaking, this transforms the original multiobjective optimization of $\f(x)$ into the single-objective set problem of finding the\del{eventual} so-called optimal $p$-distribution~\cite{auger2012hypervolume}
%\footnote{see \cite{auger2012hypervolume} for discussions on its existence.}
\begin{equation} \label{eq:optmudistribution}
		X_p^* = \argmax\limits_{\begin{array}{c} \mbox{\scriptsize $X\subseteq \R^{\dimnumber}$, \mbox{\scriptsize $|X|\leq p$}} 
		%\\[-0.5em] \mbox{\scriptsize $|X|\leq p$}
		 \end{array}} \I(X)
		% \vspace*{-0.5em}
\end{equation}
as the set of search points of cardinality $p$ (or lower) with the highest indicator value among all sets of this size \cite{auger2009theory}.

Natural candidates for practically relevant quality indicators are monotone or even strictly monotone indicators such as the epsilon-indicator \cite{zitzler2004indicator}, the R2 indicator \cite{hj1998a}, or the hypervolume indicator (\cite{zt1998b,auger2009theory}, still the only known strictly monotone indicator family to date). We remind that an indicator is called monotone if $A \preceq B \implies \I(A) \geq \I(B)$---or in other words, if it does not contradict the weak Pareto dominance relation. If $A \prec B \implies \I(A) > \I(B)$, we say that $\I$ is strictly monotone.

\paragraph{Hypervolume, Hypervolume Contribution, and Hypervolume Improvement}
Because the hypervolume indicator \cite{zt1998b,auger2009theory} and its weight\-ed variant is the only known strictly monotone indicator, we will later on use it as well in our framework. The hypervolume $HV_{\r}$ \cite{zitzler1998multiobjective}
of a finite set of solutions $S \subset \R^{\dimnumber}$ with respect to the reference point $\r\in\R^{\objnumber}$ is defined as
%\begin{equation}
%\label{hypervolume}
	$HV_{\r}(S) = \lambda_{\objnumber}\pare{ \acco{ z\in\R^{\objnumber}; \exists y \in S,\, \f(y) \prec z \prec \r } }$,
%\end{equation}
where $\lambda_{\objnumber}$ is the Lebesgue measure on the objective space $\R^{\objnumber}$ and $\f$ is the objective function.
In the case of two objective functions, the hypervolume indicator value of $p$ non-dominated solutions $S=\{s^{(1)}, \ldots, s^{(p)}\}$ with $\f_{1}(s^{(1)})\leq \f_{1}(s^{(2)})\leq\ldots\leq \f_{1}(s^{(p)})$ can also be written as the sum of the area of $p$ axis parallel rectangles: $HV_{\r}(S) = \sum_{i=1}^{p} (\f_{1}(s^{(i+1)})-\f_{1}(s^{(i)}))\cdot(r_2-\f_{2}(s^{(i)}))$; $\f_{1}(s^{(p+1)}) \egaldef r_{1}$.

Furthermore, the hypervolume {\emph contribution} $\HVC(\mathbf{s},S)$ of a search point $\mathbf{s}\in\R^{\dimnumber}$ to a solution set $S\subseteq \R^{\dimnumber}$ with respect to the reference point $\r\in\R^{\objnumber}$ is the hypervolume indicator value that we lose when we remove $\mathbf{s}$ from the set \cite{bringmann2011convergence}:
%\begin{equation}
	$\HVC_{\r}(\mathbf{s}, S) = HV_{\r}(S) - HV_{\r}(S\setminus \{\mathbf{s}\})\enspace.$
%\end{equation}

Also, the Hypervolume Improvement $\HVI_{\r}(\mathbf{s}, S)$ of a search point $\mathbf{s}\in\R^{\dimnumber}$ to a finite set $S\subset \R^{\dimnumber}$ with respect to the reference point ${\r\in\R^{\objnumber}}$ is defined as~\cite{emmerich2005emo, yang2019multi} :
%\ba\label{hypImprovement}
	$\HVI_{\r}(\mathbf{s}, S) = HV_{\r}(S\cup\acco{\mathbf{s}}) - HV_{\r}(S)\enspace.$
%\ea
Or in other words, $\HVI_{\r}(\mathbf{s}, S)$ equals the increase in hypervolume when $\mathbf{s}$ is added to the set $S$. Up to a null set $\HVI_{\r}(\mathbf{s}, S) = \HVC_{\r}(\mathbf{s}, S\cup\acco{\mathbf{s}})$.

%%%%%%%%%%%%%%%%%%%%%%%%%%%%%%%%%%%%%%%%%%%%%%%%%%%%%%%%%%%%%%%%%%%%%
\section{\SUMO:\del{ A Framework for} Building Multiobjective\del{ Algorithms} from Single-Objective \new{Algorithms}\del{Ones}}\label{sec:sofomore}
%%%%%%%%%%%%%%%%%%%%%%%%%%%%%%%%%%%%%%%%%%%%%%%%%%%%%%%%%%%%%%%%%%%%%
%As we have seen, 
Quality indicators have been introduced as a way to measure the quality of a set of objective vectors but also to define a multiobjective optimization problem as a single-objective set problem of maximizing the quality indicator as in \eqref{eq:optmudistribution}. This naturally defines a single-objective 
%With a small abuse of notation, we can also write this single-objective set problem of finding $p$ solutions with maximal indicator value as a 
$\kernumber \times \dimnumber$ dimensional problem to be maximized
\begin{equation}\label{eq:single-obj-Ind}
	%\max
	 \F: (x_1,\ldots,x_\kernumber)\in (\R^{\dimnumber})^p \mapsto \I(\{x_1, \ldots, x_\kernumber \}) \enspace.
\end{equation}
%\cheikh{Shouldn't it be $\I(\{x_1, \ldots, x_\kernumber \})$, since you defined the hypervolume from the search points space in Equation~\eqref{hypervolume} ? } Dimo: yes you are right, Cheikh, but Anne also has a point that we should come back to in the theory paper: there is a difference between indicators that can be decomposed as an indicator function in the search space and the direct mapping of all points via the objective function (like for the hypervolume here) and other indicators that cannot be decomposed (such as the inverted generational distance in search space)
%\anne{Dimo, it does not look like a too good idea to me to define the indicator with respect to points of the search space - as we always visualize it wrt to objective vectors.} Dimo: usually, an indicator is defined as mapping a set of search points to a real value, but you have a good point here (see my comment above)
Because $\dimnumber$ and in particular $\kernumber$ are typically large in practice, we usually do not attempt to solve a multiobjective optimization problem by directly optimizing \eqref{eq:single-obj-Ind}. Nevertheless, when $\I$ is the hypervolume indicator, Hern\'{a}ndez et al.\ suggest to use a Newton method to directly solve \eqref{eq:single-obj-Ind}. It assumes that $\f$ is twice continuously differentiable, in which case the gradient and Hessian of $\F$ can be computed analytically~\cite{hernandez2018set}. Yet, directly attacking \eqref{eq:single-obj-Ind} is not possible because dominated points have a zero sub-gradient and the Newton direction is therefore zero. Thus, Hern\'{a}ndez et al.\ need to start from a set of non-dominated points, close enough to the Pareto set, which requires in practice to couple the approach with another algorithm \cite{hernandez2018set}.

%But Hern\'{a}ndez et al. also realized that directly attacking \eqref{eq:single-obj-Ind} is not possible due to the fact that dominated points have a zero sub-gradient and the Newton direction is therefore zero.
Instead of directly optimizing  \eqref{eq:single-obj-Ind}, our proposed \SUMO framework performs iterative {\emph subspace optimization} of the function $\F$ and penalizes the flat landscape of $\F$ in dominated regions. 
More precisely, the basic idea behind \SUMO is to optimize $\F$ subspace- or component-wise, by iteratively fixing all but one search point $x^{(i)}$ and only optimizing the indicator with respect to $x^{(i)}$ while the other search points $X^{\neg i} := \{ x^{(1)}, \ldots, x^{(i-1)}, x^{(i+1)}, \ldots, x^{(\kernumber)} \}$ are temporarily fixed. Hence we maximize the\del{ following} functions\del{:}
\begin{equation}
	\label{eq:generalfitness}
	%\max 
	\Phi_{\I, X^{\neg i}}: x\in\R^{\dimnumber} \mapsto \I(\{x\} \cup X^{\neg i}) \enspace.
\end{equation}

If the placement of each of the $p$ search points $x^{(i)}$ ($1\leq i\leq p$) is optimized iteratively by fixing a different point set each time, as we suggest in our \SUMO framework below, we are in the setup of optimizing a dynamic fitness. More details on this aspect of our \SUMO framework will be given below in Section~\ref{sectionSofomore}.

\subsection{A Fitness Function for Subspace Optimization}
If we use as quality indicator $\I$ in \eqref{eq:generalfitness} a (strictly) monotone indicator like the hypervolume indicator, the overall fitness $\Phi$ is flat in the interior domain of regions where points are dominated.
Hence\del{ \eqref{eq:generalfitness} does not provide a reasonable fitness function allowing to find the set of $p$ solutions maximizing \eqref{eq:single-obj-Ind}.
Instead}, we suggest to not optimize \eqref{eq:generalfitness} directly but to unflatten it in dominated areas of the search space without changing the optimization goal.

\begin{figure*}~\\[-1.5ex]
\newcommand{\figwidth}{0.31}
	\centering
	\includegraphics[width=\figwidth\textwidth]{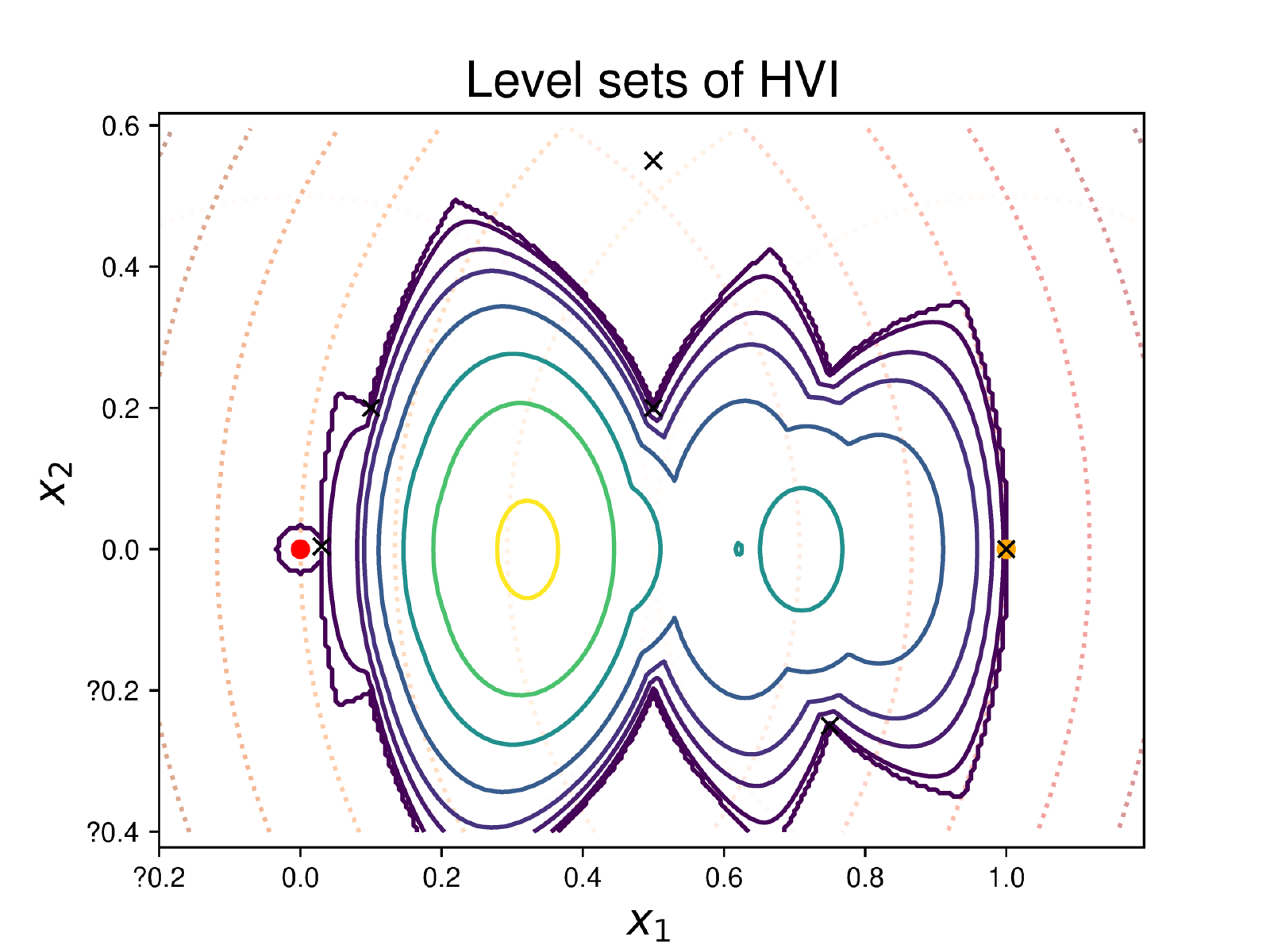}
	\includegraphics[width=\figwidth\textwidth]{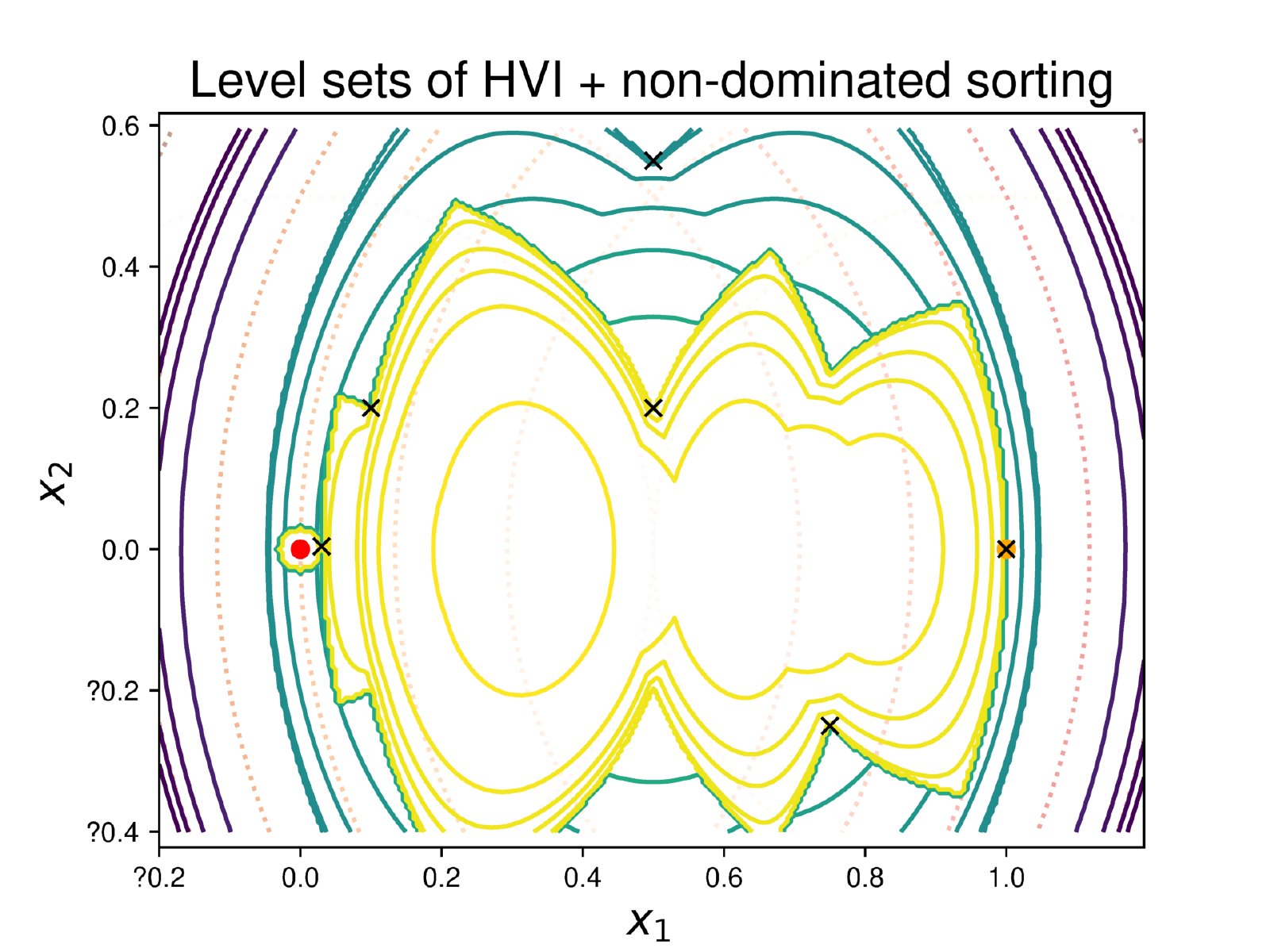}
	\includegraphics[width=\figwidth\textwidth]{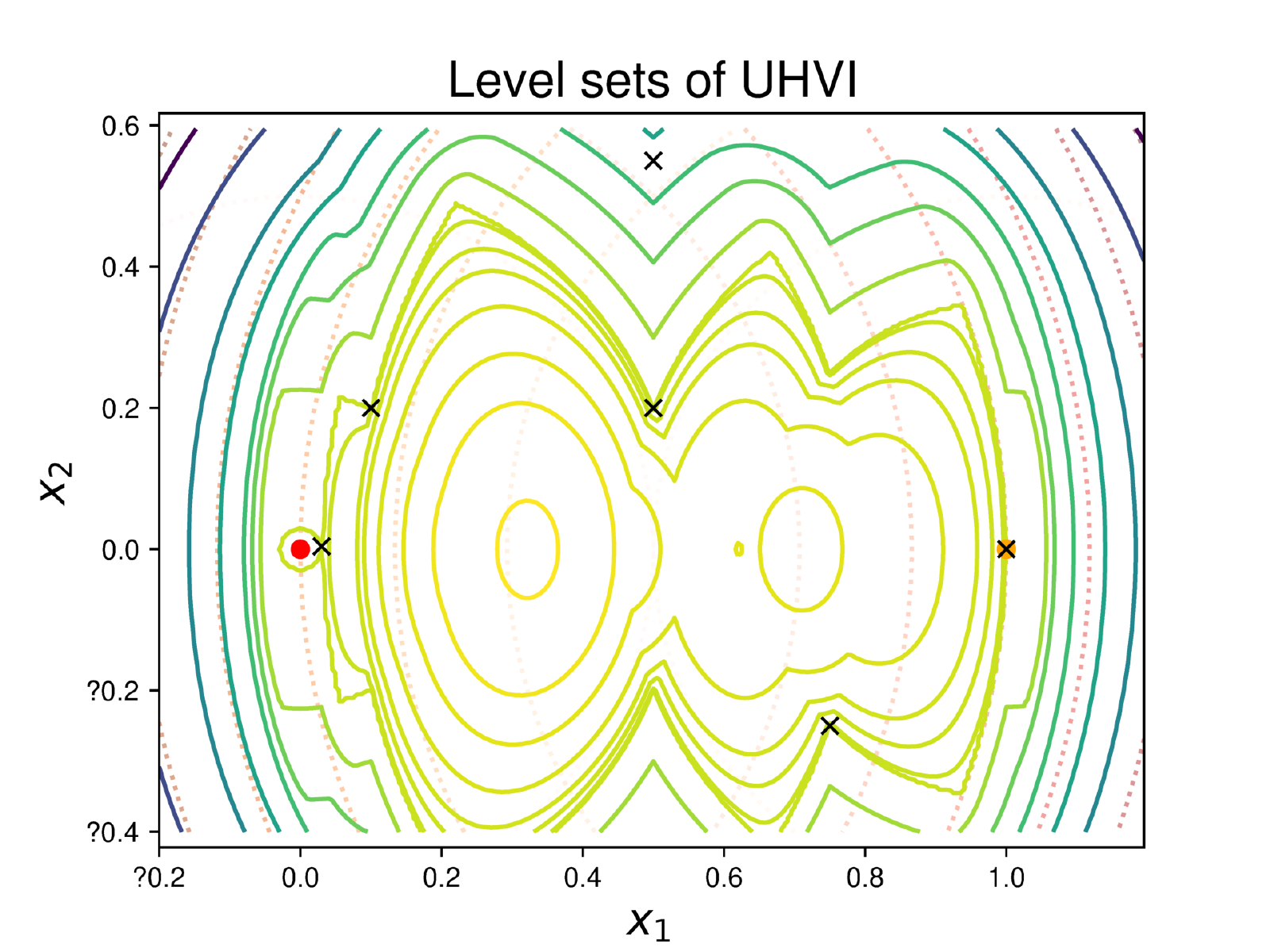}\\[-0ex]
	\includegraphics[width=\figwidth\textwidth]{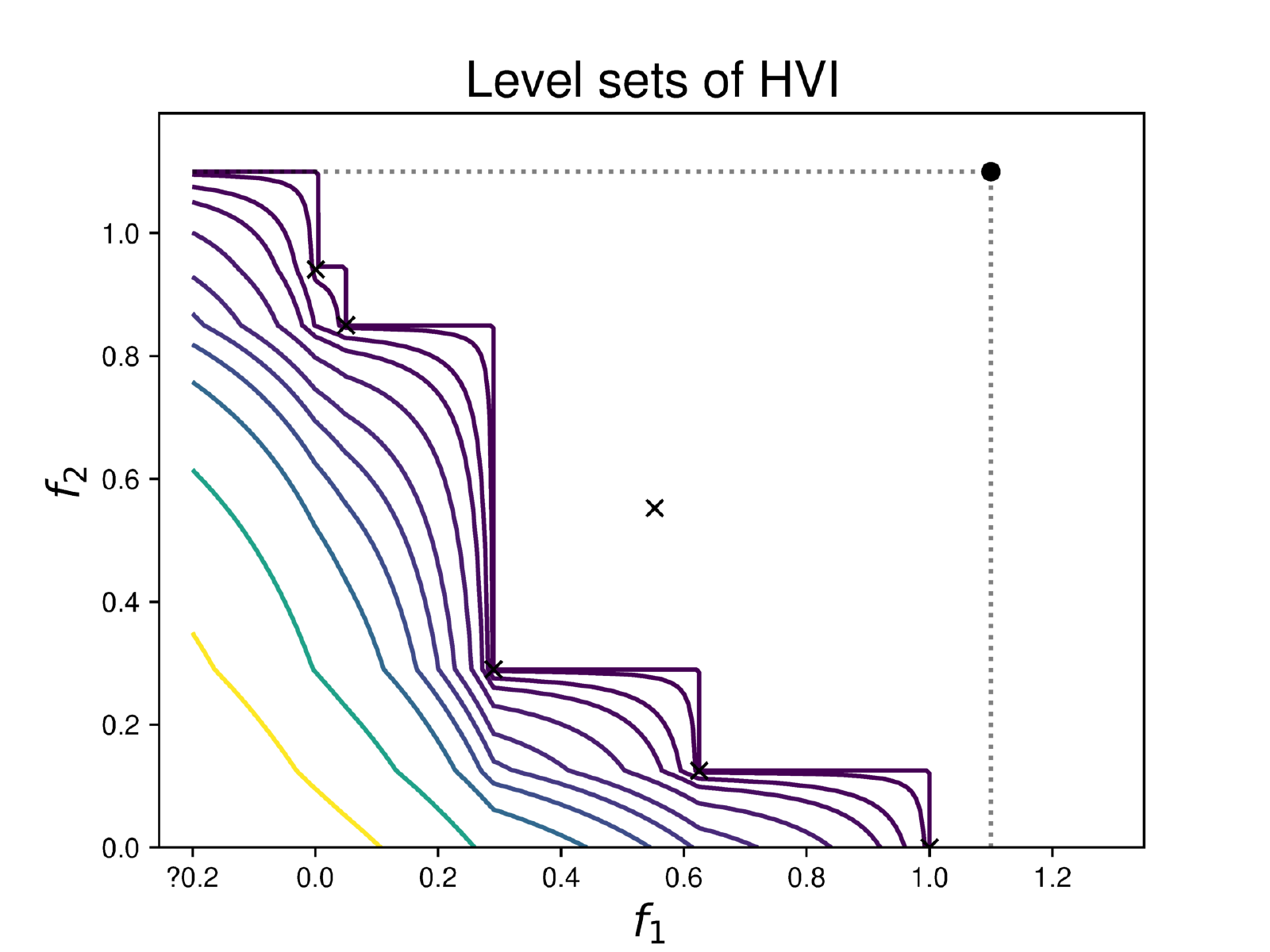}
	\includegraphics[width=\figwidth\textwidth]{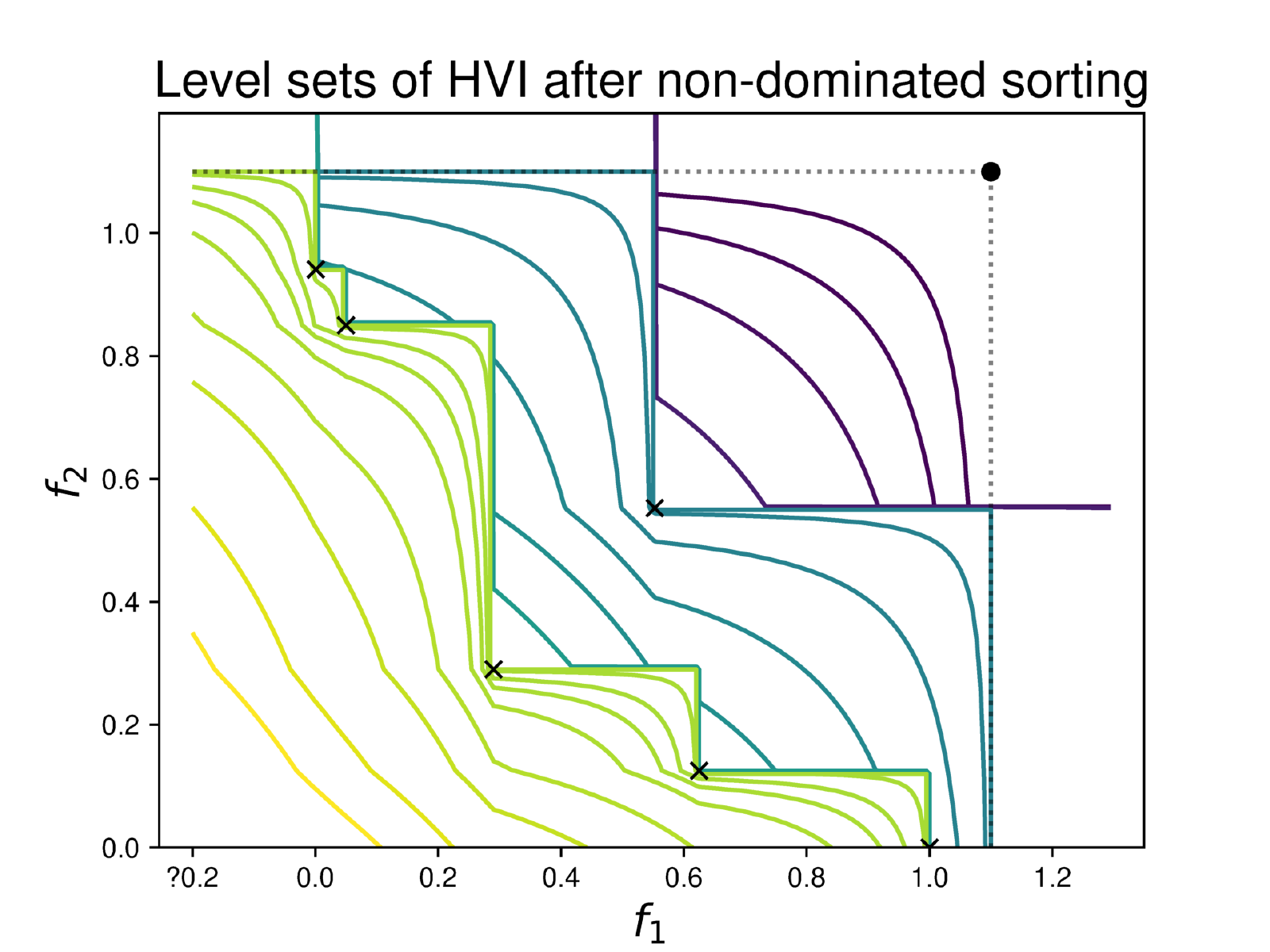}
	\includegraphics[width=\figwidth\textwidth]{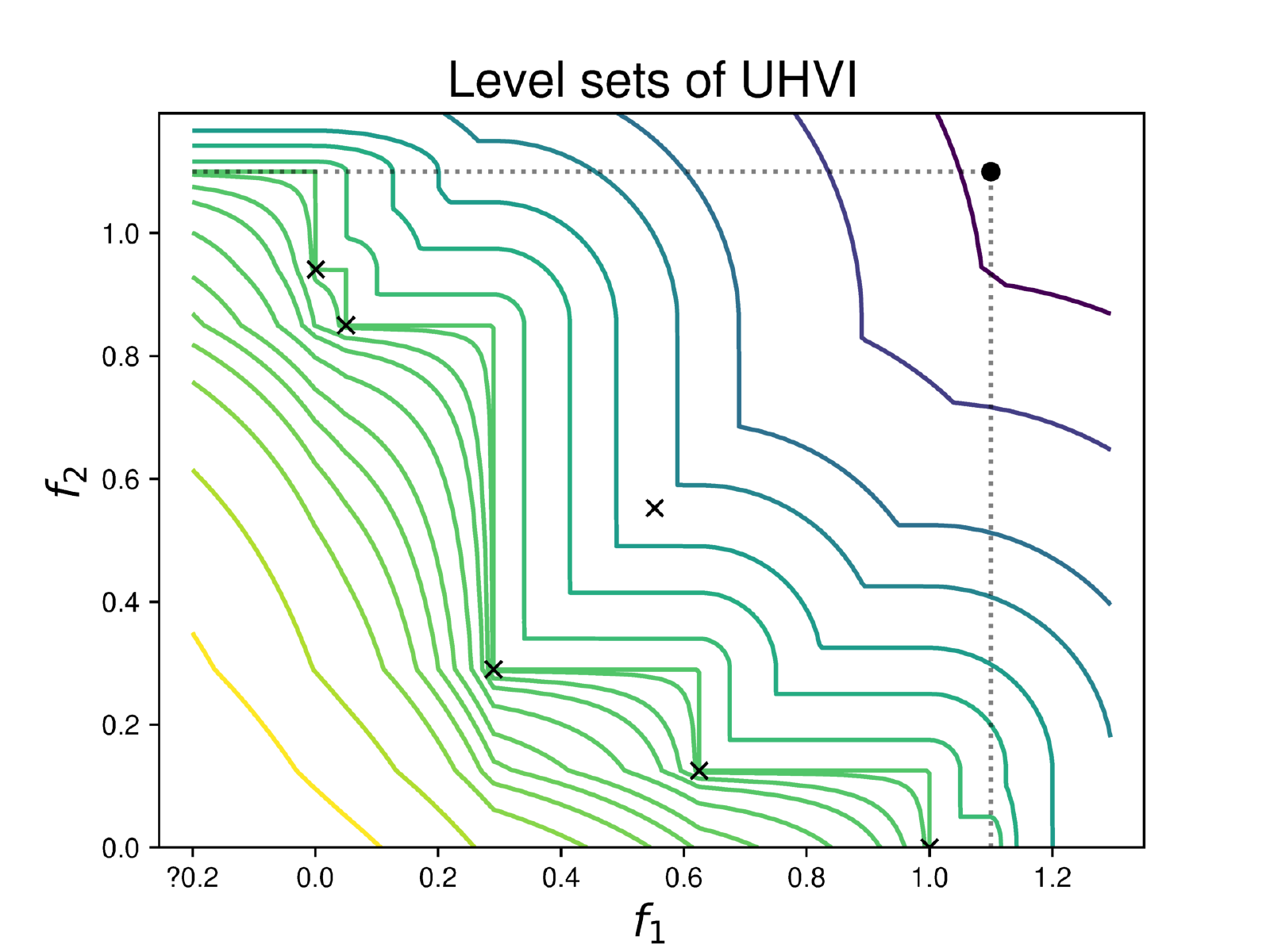}
	\vspace*{-1.5ex}
	\caption{\label{fig:fitnesscomparison}
	Comparison of block-coordinate-wise fitness functions on the double sphere problem. Given a fixed set of six solutions, $\{ [0.5, 0.2], [0.75, -0.25], [0.1, 0.2], [1,0], [0.03, 0.004], [0.5,0.55] \}$, the above three plots show level sets of equal fitness for a new search point in the search space, the second row shows the same level sets in objective space. Left: standard hypervolume improvement \HVI. Middle: \HVI within the local non-dominated fronts. Right: newly proposed hypervolume improvement (if non-dominated) together with distance to the non-dominated front (if dominated), denoted as uncrowded hypervolume improvement \UHVI. Note that the colors for the fitness levels are not comparable over indicators, but are the same for a given indicator in both search and objective space. The search space plots further show the single-objective's level sets as dotted lines. The black dot indicates the reference point of $[1.1,1.1]$.}
\end{figure*}

Any solution $x$ that is dominated by the other points in $X^{\neg i}$ will receive zero fitness $\Phi$ when we use as indicator in \eqref{eq:generalfitness} the hypervolume indicator of the entire set $\{x\}\cup X^{\neg i}$ with respect to the reference point $\r$ or replace it with the hypervolume improvement $\HVI_{\r}(x, X^{\neg i})$ of the solution $x$ to $X^{\neg i}$. This situation is depicted in the first column of Figure~\ref{fig:fitnesscomparison} where for a fixed set of six arbitrarily chosen search points, the hypervolume improvement's level sets (of equal fitness) in both search and objective space are shown. This flat fitness with zero gradient will not allow to steer the search towards better search points which has also been highlighted by Hern{\'a}ndez et al.~\cite{hernandez2018set}.

A common approach to guide an optimization algorithm in the dominated space is to use the hypervolume (contribution) as secondary fitness after non-dominated sorting~\cite{gold1989a}, as it is done for example in the SMS-EMOA~\cite{bne2007a} or the MO-CMA-ES~\cite{ihr2007a}. The idea is that all search points with a worse non-domination rank get assigned a fitness that is worse than for search points with a better non-domination rank. Within a set of the same rank, the hypervolume contribution with respect to all points with the same rank is used to refine the fitness.
The middle column of Figure~\ref{fig:fitnesscomparison} shows the resulting level sets of equal fitness. As we can see, this fitness assignment distinguishes between dominated solutions, i.e. the fitness is not flat anymore. Yet it still has another major disadvantage: the search direction in the dominated area (perpendicular to its level sets) points towards already existing non-dominated solutions. Attracting dominated solutions towards non-dominated solutions seems however undesirable, as they will compete for the same hypervolume area. Instead, we want dominated points to enter the \nnew{uncrowded}\del{ empty} space \emph{between} non-dominated points thereby complementing their hypervolume contribution (improvement).

\paragraph{Uncrowded Hypervolume Improvement}
For this purpose, we define the Uncrowded Hypervolume Improvement \UHVI based on the Hypervolume Improvement for non-dominated search points and on the Euclidean distance to the non-dominated region for dominated search points. More concretely, $\UHVI_{\r}(\mathbf{s}, S)$ of a search point $\mathbf{s}\in\R^{\dimnumber}$ with respect to a finite set $S\subset \R^{\dimnumber}$ and the reference point ${\r\in\R^{\objnumber}}$ is defined as
\begin{equation} \label{eq:UHVI}
	\UHVI_{\r}(\mathbf{s}, S) = \left\{ \begin{array}{ll} 
		\HVI_{\r}(\mathbf{s}, S) & \text{if } \emp_{S ,\r} \nprec \f(\mathbf{s})\\
		-d_{\r}(\mathbf{s}, S)   & \text{if } \emp_{S ,\r} \prec \f(\mathbf{s})
	                                                 \end{array} \right.
	    \,,
\end{equation}
where $  d_{\r}(\mathbf{s}, S) = \inf_{y\in \emp_{S,\r} } d(\f(\mathbf{s}),y)$ is the distance between an objective vector $\f(\mathbf{s})\in\R^{\objnumber}$ and the empirical non-domination front of the set $S$ defined as in \eqref{eq:empiricalnondominationfront}.

We define the fitness $\Phi_{\UHVI, X^{\neg i}}(x)$ for a search point $x\in\R^{\dimnumber}$ with respect to other solutions in $X^{\neg i}$ as
\begin{equation} \label{eq:sumofitness}
	\Phi_{\UHVI, X^{\neg i}}(x) = \UHVI_{\r}(x, X^{\neg i})\enspace.
\end{equation}
Note that $\Phi_{\UHVI, X^{\neg i}}$ is continuous on the empirical non-domination front where both the hypervolume improvement and the considered distance are zero.

\begin{figure}%
	\includegraphics[width=0.4\textwidth]{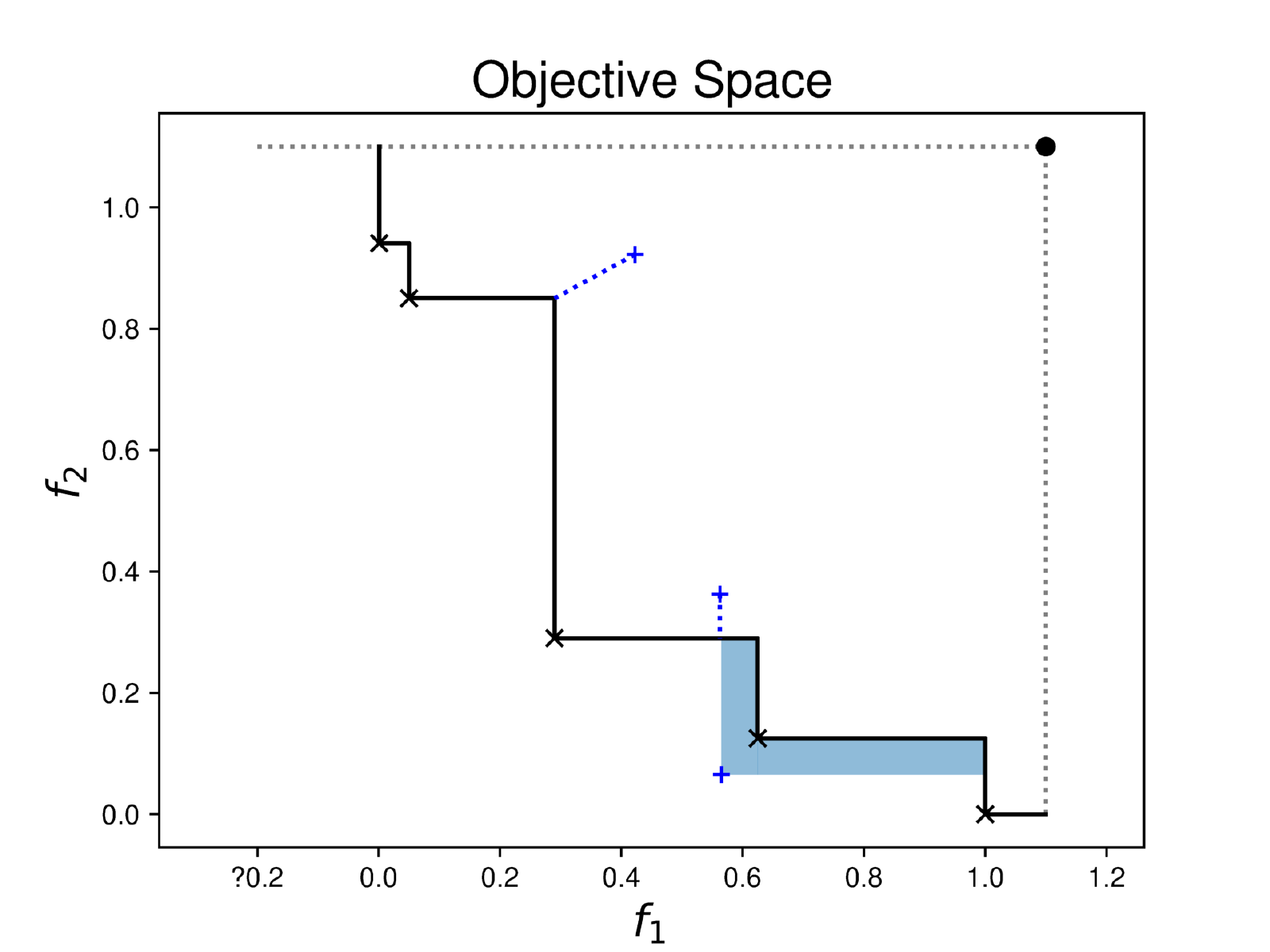}%
\vspace*{-1ex}
	\caption{\label{fig:UHVIillustration}
	Illustration of the proposed \UHVI fitness for three points in 
	objective space (blue $+$). The top two points are dominated by the given set $\mathcal{S}$ of five points (black crosses) and are thus assigned a fitness of their distance to the empirical non-dominated set while the bottom point is non-dominated and thus gets assigned the hypervolume improvement with respect to $\mathcal{S}$.}
\end{figure}

Figure~\ref{fig:UHVIillustration} illustrates this fitness for one non-dominated and two dominated search points (blue plusses) with respect to a set of six other search points (black crosses). The right-hand column of Figure~\ref{fig:fitnesscomparison} shows the level sets of this fitness.\del{ in comparison to the previously discussed fitness functions. As we can see in the example plots of Figure~\ref{fig:fitnesscomparison},} The newly introduced hypervolume improvement and distance based fitness $\Phi_{\UHVI}$ shows smooth\del{ness of its} level sets, both in search and in objective space. Maybe most importantly, in the dominated area, the fitness function's descent direction (perpendicular to its level sets) now points towards the gaps in the current Pareto front approximation.

\subsection{Iteratively Optimizing the $\Phi_{\UHVI}$ Fitness: The \SUMO Framework}\label{sectionSofomore}

After we have discussed a fitness assignment that looks worth to optimize, we come back to our initial idea of subspace optimization and define the underlying algorithmic framework behind \SUMO.

At first, we consider a \emph{single-objective} optimizer in an abstract manner as an iterative algorithm with state $\theta \in \Theta_n$ updated as $\theta_{t+1} = G^f(\theta_t,U_{t+1})$ where $f:\R^n \to \R$ is the single-objective function optimized by the optimizer and $U_{t+1}$ encodes possible random variables sampled within one iteration if we consider a randomized algorithm (and can be taken as constant in the case of a deterministic optimizer). The transition function $G^f$ contains all updates done within the algorithm in one iteration.

We assume that in each iteration $t$, the optimizer returns a best estimate of the optimum, often called incumbent solution or recommendation. This is the solution that the optimizer would return if we stop it at iteration $t$. We denote this incumbent\del{ or recommendation} as $\E: \theta \in \fonccourte{\T_{\dimnumber}}{\E(\theta) \in \R^{\dimnumber}}$---mapping the state of the algorithm to the estimate of the optimum given this state.

The overall idea behind the subspace optimization and the \SUMO framework can then be formalized as in Algorithm~\ref{alg:sumo}: after initializing $\kernumber$ single-objective algorithms with their states $\theta_{1},\dots,\theta_{\kernumber}$ and denoting their transition functions as $G_{i}^f$ ($1 \leq i \leq \kernumber$), we consider their incumbents or recommendations $\E(\theta_{i})$ as the $\kernumber$ search points that are expected to approximate the optimal $\kernumber$-distribution. 

In each step of the \SUMO framework, we choose one of the algorithms (denoted by its number $i$, with $1\leq i\leq \kernumber$) and run it $\tau_{i}$ iterations on the fitness $\Phi_{\UHVI, X^{(\neg i)}}$ to update the recommendation $x^{(i)}$ while keeping all other recommendations fixed. It is important to note that the fitness used for algorithm $i$ is actually changing \emph{ dynamically} with the optimization because it depends on all the other incumbents but $x^{(i)}$ which, over time, are expected to move towards the Pareto set as well.

Algorithm~\ref{alg:sumo} proposes a generic framework where the order in which the single-objective algorithms are run and the number of iterations for them are not explicitly defined. A simple strategy would be to choose the algorithms at random or in a given, fixed order and run each single-objective algorithm a fixed number of time steps. But also more elaborate strategies can be envisioned, for example based on the idea of multi-armed bandits \cite{bc2012a}: we can log the changes in the fitness value of each incumbent over time and favor as the next chosen algorithms the ones that give the highest expected fitness improvements. Note also that the single-objective algorithms' types may be different such that we can combine local with global algorithms or even change the algorithms over time, allow restarts etc. In the following experimental validation of our concept, however, we choose a single optimization algorithm and a simple, random strategy to choose which of them to run next.
\del{Finally, not even the usage of the proposed hypervolume and distance based fitness is fixed; also other fitness functions, based on a monotone quality indicator and a penalty for dominated points can be envisioned.}
\begin{algorithm}
\caption{General \SUMO Framework, with fitness of \eqref{eq:sumofitness}}
\label{alg:sumo}
\begin{algorithmic}[1]
	\State {\bf Given: } the initial states of $\kernumber$ single-objective optimizers, $\theta_{1},\dots,\theta_{\kernumber}$
	\State Initialize incumbents: $X = \{x^{(1)}=\E(\theta_{1}), \ldots, x^{(p)}=\E(\theta_{p})\}$
	\While{\textbf{not} stopping criterion met}
		\State Choose $i$ in $\acco{1,\dots,p}$ and $\tau_{i}\in\N$
		% \State Choose an integer $\tau_{i}$
		\State {\bf REPEAT} $\tau_{i}$ times:
		\State \quad $\theta_{i} \leftarrow G_{i}^{\Phi_{\UHVI, X^{(\neg i)}}}(\theta_{i},U_{i})$ 
		\Comment{run \ensuremath{i}th algo on fitness
		    \ensuremath{\Phi_{\UHVI, X^{(\neg i)}}}}
		\State $x^{(i)} \leftarrow \E(\theta_{i})$ \Comment{update $x^{(i)}$ in $X$}
	\EndWhile
	\State \Return $x^{(1)},\dots,x^{(\kernumber)}$
\end{algorithmic}
\end{algorithm}

With a simple change, Algorithm~\ref{alg:sumo} can be made parallelizable (resulting in slightly different search dynamics though): postponing the updates of the $x^{(i)}$ after every algorithm has been touched at least once makes the optimization of the fitness functions independent such that they can be performed in parallel.

\paragraph{Relation of \SUMO with other existing algorithms}
We briefly discuss how some existing algorithms and algorithm frameworks relate to the new \SUMO proposal.

The coupling of single-objective algorithms to form a multiobjective one has been done before, especially in the MOEA/D framework \cite{zl2007a}. In MOEA/D, $p$ static search directions (in objective space) are defined via $p$ (single-objective) scalarizing functions. Each of them is optimized in parallel with solutions potentially shared between neighboring search directions.
On the contrary, the fitness in \SUMO is dynamic, depending on the other incumbents. Optimizing a set of scalarizing functions in classical approaches to multiobjective optimization have static optimization problems to solve without any interaction between them \cite{miet1999a}.

Many other EMO algorithms, such as NSGA-II, SMS-EMOA, or MO-CMA-ES are not covered by the \SUMO framework. One simple reason is that the \UHVI is newly defined.

The already mentioned Newton algorithm on the hypervolume indicator fitness of \cite{hernandez2018set} is probably the closest existing approach from \SUMO, but \cite{hernandez2018set} needs to initialize the Newton algorithm with a set of non-dominated solutions in order for the algorithm to optimize due to the flat regions of its objective space. Also algorithms for expensive multiobjective optimization based on the optimization of the expected hypervolume improvement \cite{wedp2010a} can be seen as related to \SUMO, although the proposal of new solutions in algorithms like SMS-EGO \cite{pwbv2008a} or S-metric based ExI \cite{ek2008a} use Gaussian Processes to model the objective function.
These algorithms, in contrary to \SUMO, propose iteratively a single solution based on the \emph{expected} hypervolume improvement over all known solutions and do not aim at replacing successively a single recommendation by another (better) one. Interesting to note is that algorithms like SMS-EGO and S-metric based ExI
employ the expected hypervolume indicator improvement as fitness while the approach of Keane \cite{kean2006a} ``uses the Euclidean distance to the nearest vector in the Pareto front'' \cite{wedp2010a}.

\section{\COMOCMA\del{ Instantiating \SUMO\ with hypervolume, distance penalization, and the CMA-ES}}\label{sec:comocma}
%%%%%%%%%%%%%%%%%%%%%%%%%%%%%%%%%%%%%%%%%%%%%%%%%%%%%%%%%%%%%%%%%%%%%
%
In this section, we instantiate \SUMO\ with the CMA-ES as single objective optimizer.

Regarding the choice of which optimization algorithm to run (and how long), we opt for a simple strategy: 
we sample a permutation from $S_p$, the set of all permutations on $\{1,\ldots,p\}$ uniformly at random and use this fixed permutation to touch each algorithm $i$ once in the order of the permutation. Once all algorithms have been touched, we then resample a new permutation.
%Each $p$ iterations\anne{not so clear what should be an iteration for the MO algorithm, maybe looping over p \kernels - in any case this need to be clarified}, we sample a new permutation from $S_p$, the set of all permutations on $\{1,\ldots,p\}$ uniformly at random and use this fixed permutation to touch each algorithm $i$ once in the order of the permutation. 
We run each algorithm for a single iteration. Letting the algorithms run for a too long period right from the start seems suboptimal. As the fitness is dynamic, we do not need to optimize it too precisely.
We mainly have two requirements for the choice of single objective optimizers: (i) an optimization algorithm has to be stoppable at any iteration and resumable thereafter and (ii) an optimization algorithm needs to be able to give a good recommendation about the best estimate of the optimum, given its current state. The Covariance Matrix Adaptation Evolution Strategy (CMA-ES, \cite{ho2001a}) is a natural choice. Not only is it a state-of-the-art algorithm for difficult blackbox optimization problems but also does it fulfill our requirements. In CMA-ES, the state of the algorithm $\theta$ is composed of a step size ${\bf \sigma}$ and the parameters of a multivariate normal distribution, namely a mean vector ${\bf m}$ representing the favorite solution and a covariance matrix ${\bf C}$.
In addition, two $n$-dimensional evolution paths speed up step-size and covariance matrix adaptation. For each $\theta_i$, the incumbent solution $x^{(i)} = \E(\theta_i)$ is the mean ${\bf m}$ of the CMA algorithm.
A convenient implementation of CMA-ES is via the ask and tell interface \cite{collette2010object}, where the ask function returns\del{ given the state $\theta$ of the algorithm} $\lambda$ candidate solutions\del{ sampled and to be evaluated} and the tell function updates the state from their fitness values.
The interface allows to easily stop and resume the optimization and to integrate the dynamic fitness of \SUMO, see Algorithm~\ref{alg:como}.
We call this instantiation of the \SUMO framework \COMOCMA. The $p$ CMA-ES instances are called \kernels.

\begin{algorithm}
\caption{The \COMOCMA: an instance of the \SUMO framework with the CMA-ES as single-objective optimizer}
\label{alg:como}
\begin{algorithmic}[1]
	\State {\bf Required:}
	\State \quad objective function $\f = (\f_{1},\ldots, \f_{\objnumber})$ in dimension $\dimnumber$
	\State \quad lower and upper bounds for each variable $\text{\ttfamily lower}\in\R^{\dimnumber}$
	\Statex\quad\quad and $\text{\ttfamily upper}\in\R^{\dimnumber}$ of a region of interest
  \State \quad number of desired solutions $\kernumber$
	\State \quad global initial step-size $\sigma_{0}$ for all CMA-ES
	\State \quad fixed reference point $\r$ for the hypervolume indicator\vspace{0.5em}
	\State \textbf{Initialization:}
		\State $x^{(i)} = \text{uniformSample}(\text{\ttfamily lower},\text{\ttfamily upper})$ for all $1\leq i \leq \kernumber$
		\State evaluate all $x^{(i)}$ on $\f$ and store the $\f(x^{(i)})$ for later use
		\State $\text{es}^{(i)} \leftarrow \pare{\mu/\mu,\lambda}\text{-CMA-ES}\pare{x^{(i)},\sigma_{0}}$ for all $1\leq i \leq \kernumber$\vspace{0.5em}
	
		\While{\textbf{not} stopping criterion}

			\State sample uniformly at random a permutation $\pi$ from all 
			\Statex\quad\quad\quad permutations on $\{1,\ldots,\kernumber\}$
		
			\For{i = 1 to \kernumber}
				\State $\{s^1, \ldots, s^\lambda\} \leftarrow \text{es}^{(\pi(i))}\text{.ask()}$ \Comment{get $\lambda$ offspring from}
				\Statex\Comment{$\pi(i)$th CMA-ES}
				\State compute the fitness $\Phi(s^{j}) = \Phi_{\UHVI, X^{(\neg \pi(i))}}(s^{j})$ 
				\Statex\quad\quad\quad\quad for all $1\leq j\leq \lambda$ 
\del{				\State sort the $\{s^1, \ldots, s^\lambda\}$ according to $\Phi$ such 
				\Statex\quad\quad\quad\quad that $\Phi(s^{1:\lambda})\leq \ldots\leq \Phi(s^{\lambda:\lambda})$
				\State $\text{es}^{(\pi(i))}.\text{tell}(s^{1:\lambda}, \ldots, s^{\lambda:\lambda}, \Phi(s^{1:\lambda}), \ldots, \Phi(s^{\lambda:\lambda}))$
}
				% \Comment{update \text{es}}
				\State $\text{es}^{(\pi(i))}.\text{tell}(\del{s^{1}, \ldots, s^{\lambda}, }\Phi(s^{1}), \ldots, \Phi(s^{\lambda}))$
				\State $x^{(\pi(i))} \leftarrow \text{es}^{(\pi(i))}.\text{mean}$
				\State update the stored objective vector $\f(x^{(\pi(i))})$\label{alg2line18}
\del{				\Statex\quad\quad\quad\quad with $\f(x^{(\pi(i))})$}

		\EndFor
		\EndWhile
		\State {\bf Return} $x^{(1)},\dots, x^{(\kernumber)}$

\end{algorithmic}
\end{algorithm}

We see in particular how \del{the single-objective }CMA-ES is integrated into \SUMO via its ask-and-tell interface. After choosing the next \kernel $i$, the corresponding CMA-ES instance samples $\lambda$ solutions (``ask''). It then evaluates them on the uncrowded hypervolume improvement based fitness defined in Eq~\eqref{eq:sumofitness}---given all other \kernels being fixed. After sorting the $\lambda$ solutions with respect to their fitness, \COMOCMA feeds the sampled points with their fitness values back to the CMA-ES instance (``tell'') which updates all its internal algorithm parameters. Finally, the new mean of the corresponding CMA-ES instance updates the list of the \COMOCMA's proposed $\kernumber$ solutions. Note here that CMA-ES is usually not evaluating the mean of the sample distribution which therefore is done in line~\ref{alg2line18}.
\section{Experimental Validation}
\label{sec:experiments}
\newcommand{\pycma}{\href{https://github.com/CMA-ES/pycma}{\texttt{pycma}}}
We present in this section numerical experiments of the \COMOCMA. Though, in principle, the algorithm can be defined for any number of $\objnumber$ objectives, we present results only for $\objnumber=2$\del{ here}.
We use the \pycma\ Python package \cite{hansen2019pycma} version $2.6.0$ for CMA-ES as single-objective optimizer without further parameter tuning.

\subsection{Test Functions and Performance Measures}
For a matrix $P$ and two vectors $x$ and $y$, we denote 
\ba\label{quad}
\textbf{Quad}(P,x,y) = \pare{x-y}^\top P \pare{x-y}. 
\ea
We also denote by $\mathbf{0}$ the all-zeros vector, $\mathds{1}$ the all-ones vector, and $e_{k}$ the unit vector with its only nonzero value\del{ (``$1$'') being} at position $k$.
 Starting from a positive diagonal matrix $\dsp\Delta$, and two independent orthogonal matrices
$O_{1}$ and $O_{2}$, we consider the classes of bi-objective convex quadratic problems \textbf{Sep-$k$}, \textbf{One} and \textbf{Two} defined as follows ~\cite{toure2019bi}
\begin{itemize}
\item   $\f_{1,\Delta}^{ \text{\textbf{sep-$k$}}}(x) =  \frac{ \textbf{Quad}(\Delta,x,\mathbf{0}) }{ \textbf{Quad}(\Delta,\mathbf{0},e_{k}) }$, $\f_{2,\Delta}^{ \text{\textbf{sep-$k$}}}(x) =  \frac{ \textbf{Quad}(\Delta,x,e_{k}) }{ \textbf{Quad}(\Delta,\mathbf{0},e_{k}) }$.
\item $ \f_{1,\Delta}^{\text{\textbf{one}}}(x) =  \frac{ \textbf{Quad}(O_{1}^{T}\Delta O_{1},x,\mathbf{0}) }{ \textbf{Quad}(O_{1}^{T}\Delta O_{1},\mathbf{0},\mathds{1}) }$, $\f_{2,\Delta}^{\text{\textbf{two}}}(x)=  \frac{ \textbf{Quad}(O_{1}^{T}\Delta O_{1},x,\mathds{1}) }{ \textbf{Quad}(O_{1}^{T}\Delta O_{1},\mathbf{0},\mathds{1}) }$
\item $ \f_{1,\Delta}^{\text{\textbf{two}}}(x) =  \frac{ \textbf{Quad}(O_{1}^{T}\Delta O_{1},x,\mathbf{0}) }{ \alpha }$, $ \f_{2,\Delta}^{\text{\textbf{two}}}(x) =  \frac{ \textbf{Quad}(O_{2}^{T}\Delta O_{2},x,\mathds{1}) }{ \alpha }$
with $\alpha = \max\pare{ \textbf{Quad}(O_{1}^{T}\Delta O_{1},\mathbf{0},\mathds{1}), \textbf{Quad}(O_{2}^{T}\Delta O_{2},\mathbf{0},\mathds{1}) }$.
\end{itemize}
If $\Delta$ is the identity matrix, we call the problems as \textbf{sphere-sep-$k$} in the first case and $\textbf{bi-sphere}$ in the second and third cases (the rotations are ineffective).
If $\Delta(i,i) = 10^{6\frac{i-1}{\dimnumber-1}}$ for $i=1,\dots,\dimnumber$, then we denote the problems as \textbf{elli-sep-$k$}, \textbf{elli-one} or \textbf{elli-two}. If $\Delta(1,1) = 10^{-4}$, $\Delta(2,2) = 10^{4}$ and $\Delta(i,i) = 1$ for $i=3,\dots,n$, then we have \textbf{cigtab-sep-$k$}, \textbf{cigtab-one} or \textbf{cigtab-two}.

We fix the reference point to $\r = \pare{1.1, 1.1}$. The scalings above ensure that the reference point is dominated by all Pareto fronts considered, and that the \textbf{Sep-$k$} and the \textbf{One} problems have the same Pareto front (see ~\cite{toure2019bi}) than the bi-sphere $g: \foncfast{\croc{0,1}}{t}{(1-\sqrt{t})^{2}}{\R}$. Note that the expression does not depend on the dimension $\dimnumber$. 

\begin{figure*}
	\centering
	\includegraphics[width=0.40\columnwidth]{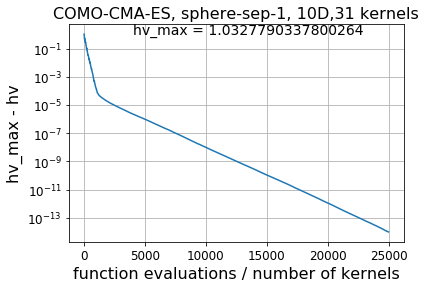}
		\hspace{-0.2cm}
		\includegraphics[width=0.40\columnwidth]{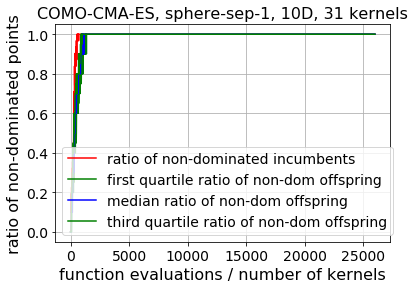}
		\hspace{-0.2cm}
	\includegraphics[width=0.40\columnwidth]{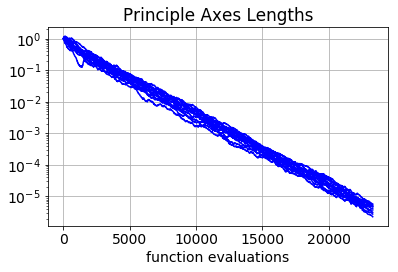}
	\hspace{-0.2cm}
		\includegraphics[width=0.40\columnwidth]{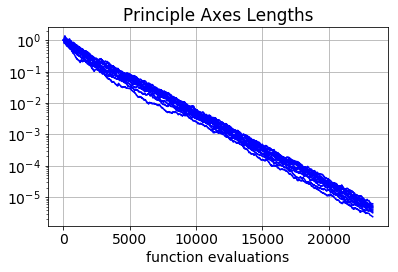}
\hspace{-0.2cm}
		\includegraphics[width=0.40\columnwidth]{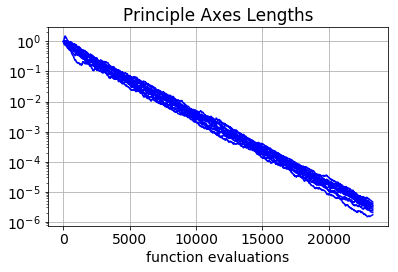}
		\hspace{-0.2cm}
			\includegraphics[width=0.40\columnwidth]{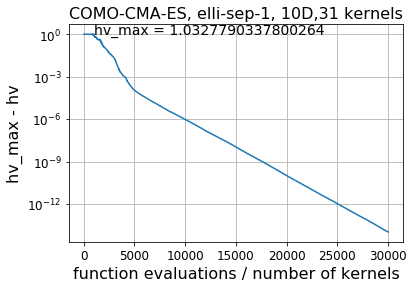}
		\hspace{-0.2cm}
	\includegraphics[width=0.40\columnwidth]{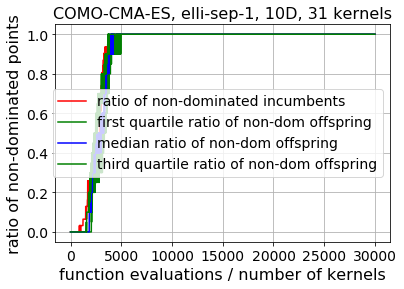}
\hspace{-0.2cm}	
	\includegraphics[width=0.40\columnwidth]{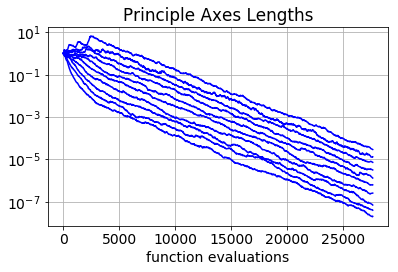}
	\hspace{-0.2cm}
		\includegraphics[width=0.40\columnwidth]{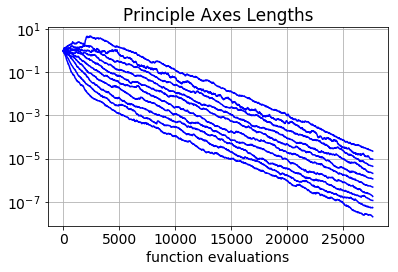}
\hspace{-0.2cm}
		\includegraphics[width=0.40\columnwidth]{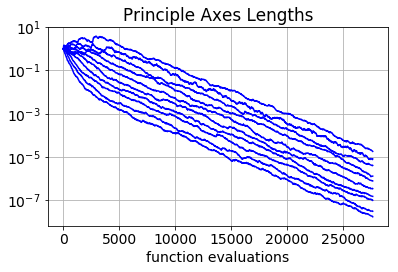}
\hspace{-0.2cm}	
	\includegraphics[width=0.40\columnwidth]{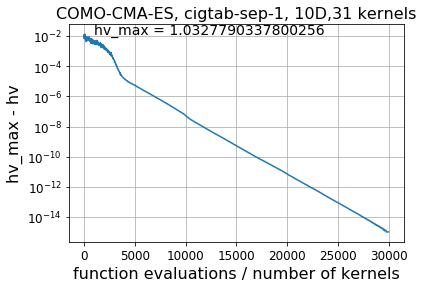}		
\hspace{-0.2cm}	
	\includegraphics[width=0.40\columnwidth]{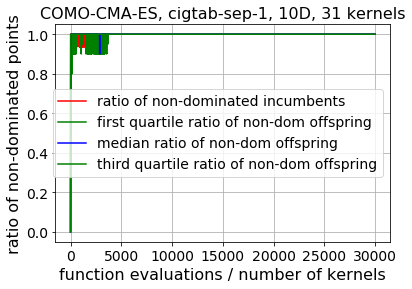}
\hspace{-0.2cm}	
	\includegraphics[width=0.40\columnwidth]{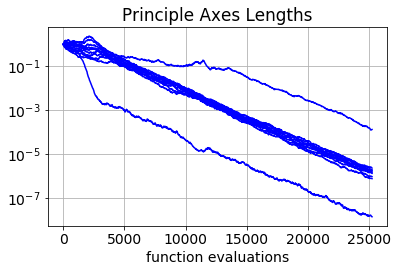}
	\hspace{-0.2cm}
		\includegraphics[width=0.40\columnwidth]{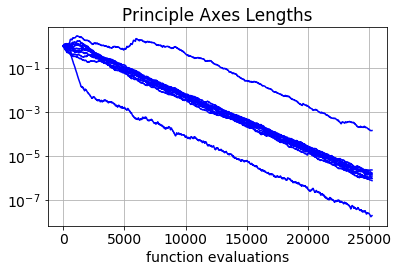}
\hspace{-0.2cm}
		\includegraphics[width=0.40\columnwidth]{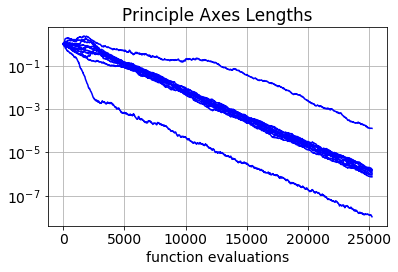}
		\hspace{-0.2cm}	
	\includegraphics[width=0.40\columnwidth]{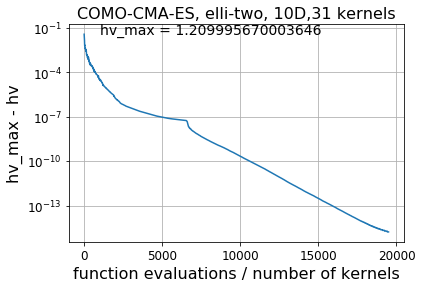}		
\hspace{-0.2cm}	
	\includegraphics[width=0.40\columnwidth]{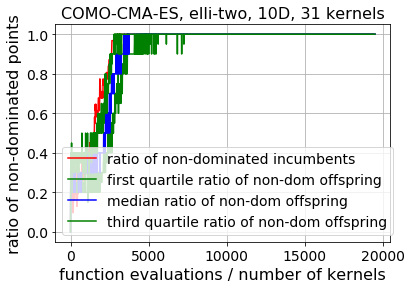}
\hspace{-0.2cm}	
	\includegraphics[width=0.40\columnwidth]{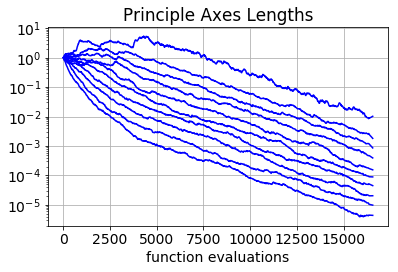}
	\hspace{-0.2cm}
		\includegraphics[width=0.40\columnwidth]{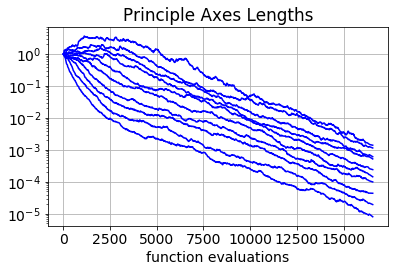}
\hspace{-0.2cm}
		\includegraphics[width=0.40\columnwidth]{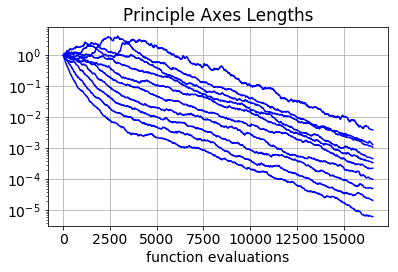}
	\caption{\label{fig:comocmaes} 
%	Convergence of the COMO-CMA-ES on \textbf{sphere-sep-$1$} (first row), \textbf{elli-sep-$1$} (second row), \textbf{cigtab-sep-$1$} (third row) and \textbf{elli-two} in $10$D with $31$ \kernels. Each row represents respectively the \convergencegap, the ratio of non-dominated points and the $3$ square root eigenspectra of the \kernels' covariance matrices.
Convergence of COMO-CMA-ES on \textbf{sphere-sep-$1$} (first row), \textbf{elli-sep-$1$} (second row), \textbf{cigtab-sep-$1$} (third row) and \textbf{elli-two} in $10$D with $31$ \kernels. The first column represents the \convergencegap. The second column is the ratio of non-dominated points among the \kernels incumbents (red) and the quartiles of the $31$ ratios of non-dominated points among each \kernel's incumbent and its offspring (median in blue and the remaining quartiles in green). And the last three columns are square root eigenspectra of uniform randomly chosen $3$ among the $31$ \kernels' covariance matrices.
}
\end{figure*}

 We use two performance measurements in each run of an algorithm. First the
 %\begin{itemize}
% \item 
 {\bf \convergencegap} defined as  the difference \del{in $\log$-scale }between an offset called \texttt{hv\_max} and the hypervolume of the $\kernumber$ points $\{x^{(1)}, \ldots,x^{(\kernumber)}\}$ found by the algorithm (in case of \COMOCMA or of the population for the other algorithms tested) called \texttt{hv}; and second the
 %\item 
{\bf \archivegap} defined as the difference \del{in $\log$-scale }between an offset called \texttt{hvarchive\_max} and the hypervolume of all non-dominated points found by the algorithm called \texttt{hvarchive}.
% \end{itemize}
The setting of \texttt{hv\_max} is done for each problem as the maximum hypervolume value of $\kernumber$ \kernels found so far anytime the problem was optimized in our machines, plus a small number ($< 10^{-14}$). 
For the \textbf{Sep-$k$} and the \textbf{One} problems, we take \texttt{hvarchive\_max} as $1.1^{2} - \int_{0}^{1}g(t)\diff t = 1.21 - \frac16$ which corresponds to the hypervolume of the theoretical Pareto front. For the two-class of problems, we use the analytic expression of their Pareto set ~\cite{toure2019bi} to sample a large number of points on the Pareto set, and compute their hypervolume as \texttt{hvarchive\_max}. Thus for the \textbf{elli-two} problem in dimension $10$, we sample $10^{7}$ points.

\subsection{Linear convergence of \COMOCMA}
We investigate the convergence of \COMOCMA for different dimensions and number of kernels, and display the results on the \textbf{sphere-sep-$1$}, \textbf{elli-sep-$1$}, \textbf{cigtab-sep-$1$} and \textbf{elli-two} functions for $\dimnumber=10$ and $\kernumber=31$. The global initial step-size is set to $\sqrt{10}$ and the \nnew{initial} lower, upper bounds (line~8 of Algorithm~\ref{alg:como}) respectively to $-5\mathds{1}$, $5\mathds{1}$.
In Figure~\ref{fig:comocmaes}, we observe linear convergence in the \convergencegap (first column) on all test functions, starting roughly when all displayed ratios of non-dominated points reach $1$ (second column). The last three columns of Figure~\ref{fig:comocmaes} illustrate the eigenspectra of the \kernels covariance matrices.
The first two columns reveal two phases.

First, the \kernels incumbents approach the non-dominated region: for \textbf{sphere-sep-$1$} this takes about $1500$ evaluations per \kernel, for \textbf{elli-sep-$1$}, \textbf{cigtab-sep-$1$} and $\textbf{elli-two}$ it takes about $5000$, $4000$ and $6600$ evaluations per \kernel.
Afterwards, the \convergencegap converges linearly. In our settings, there are $11$ evaluations per \kernel during the update of a \kernel, thus for the \textbf{*-$1$} functions (which have the same Pareto set and front), the linear convergence rate is about $10^{-\frac{6\times 11}{15000}}$ and for $\textbf{elli-two}$, it is about $10^{-\frac{3\times 11}{5000}}$.

For the first $1000$ function evaluations per \kernel on \textbf{elli-sep-$1$}, there is no point dominating the reference point, which means that the algorithm started far from the Pareto front.
Looking at \textbf{elli-two}, we confirm that it has a different Pareto front than the three other problems: $\texttt{hv\_max} = 1.2099...$ (instead of $1.0327...$).

The Uncrowded Hypervolume Improvement depends on other \kernels' incumbents and therefore changes in each iteration. Yet, the last three columns are similar to what one would observe when optimizing a single objective convex-quadratic function with corresponding Hessian matrix.
After a large enough number of iterations, the probability that the incumbents and their offspring are in the Pareto set becomes close to $1$.
Then if the incumbents $X = \acco{x^{(1)}, \ldots x^{(\kernumber)} }$ is a subset of the Pareto set and $x$ is non-dominated,  Eq~\eqref{eq:sumofitness} becomes: $\Phi_{\UHVI, X^{\neg i}}(x) = \UHVI_{\r}(x, X^{\neg i}) = \HVI_{\r}(x, X^{\neg i})$.
Its Hessian on smooth bi-objective problems is
$\nabla^{2} \Phi_{\UHVI, X^{\neg i}}(x) = \pare{\f_{2}(x)-\f_{2}( x^{(i-1)})}\nabla^{2} \f_{1}(x) +
\pare{\f_{1}(x)-\f_{1}( x^{(i+1)} ) }\nabla^{2} \f_{2}(x) +{}$\linebreak[4]
$\nabla \f_{2}(x)\nabla \f_{1}(x)^\top + \nabla \f_{1}(x)\nabla \f_{2}(x)^\top$.
For our test functions, it is a linear combination of the single objectives Hessian matrices, up to a rank-one matrix and its transpose (the gradients are colinear on the Pareto set of bi-objective convex quadratic problems~\cite{toure2019bi}). That might give a glimpse on the behaviour seen in the last three columns of Figure~\ref{fig:comocmaes}.

\begin{figure*}
	\centering
	\includegraphics[width=0.24\textwidth]{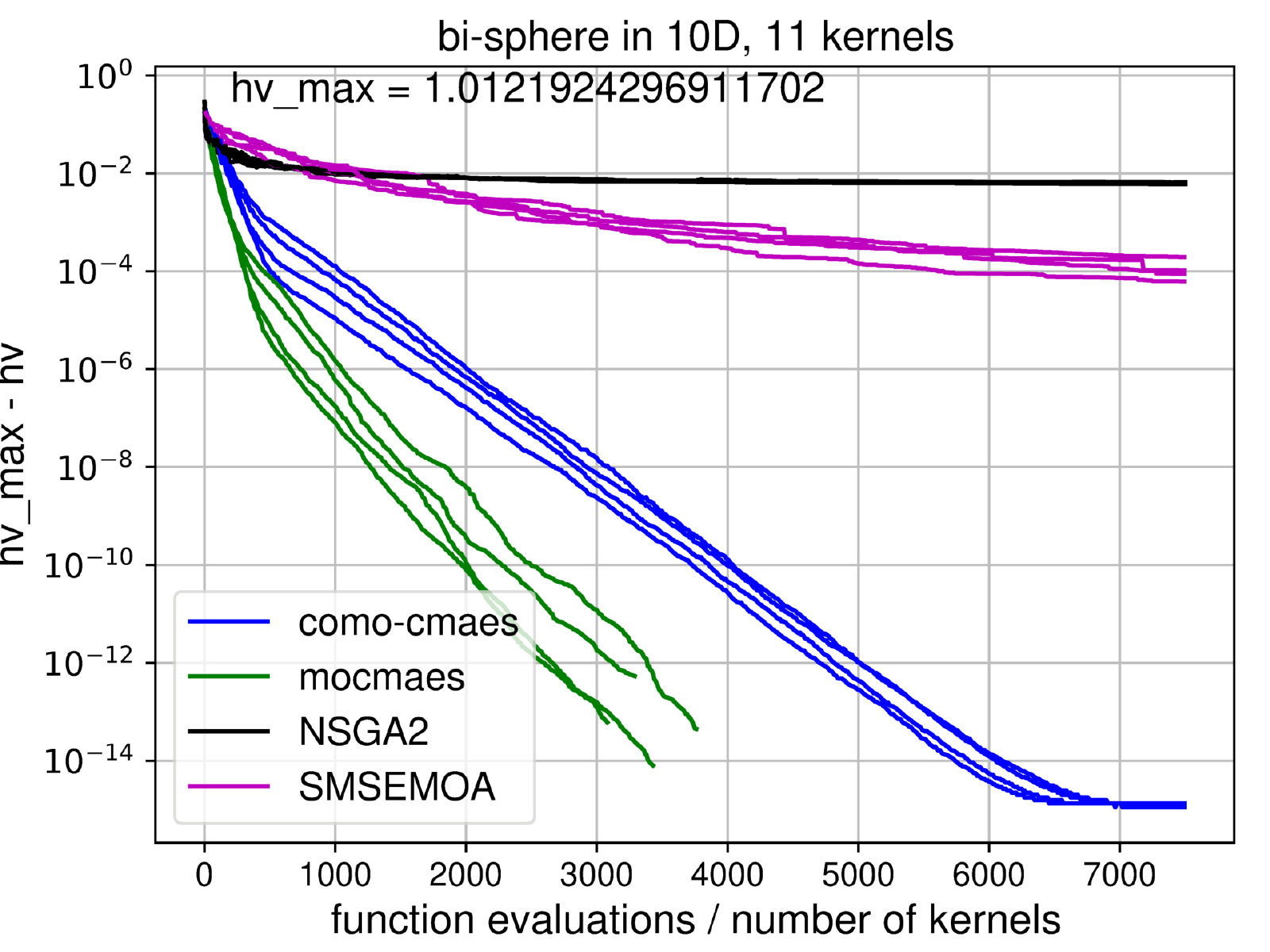}
	\includegraphics[width=0.24\textwidth]{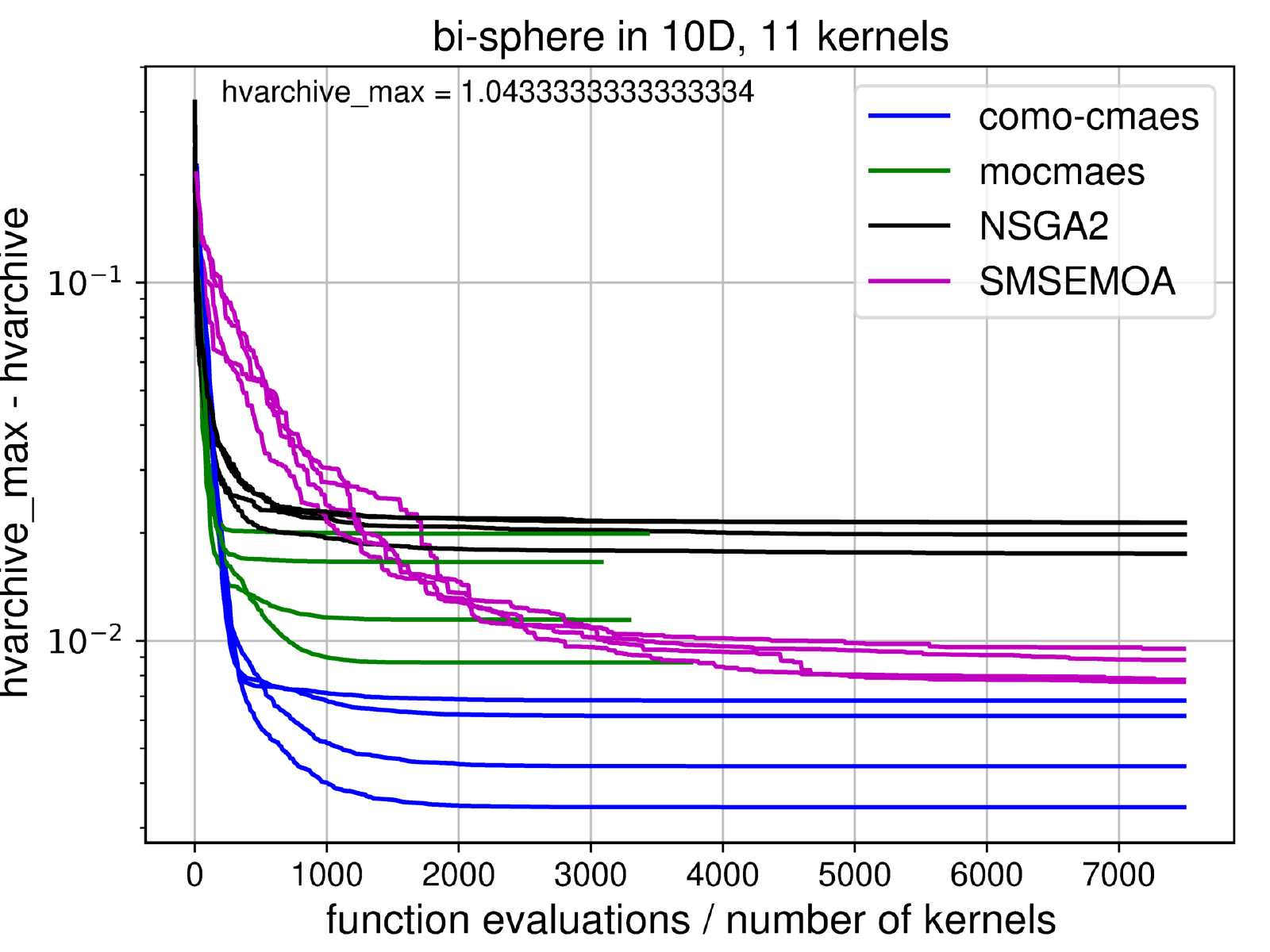}
	\includegraphics[width=0.24\textwidth]{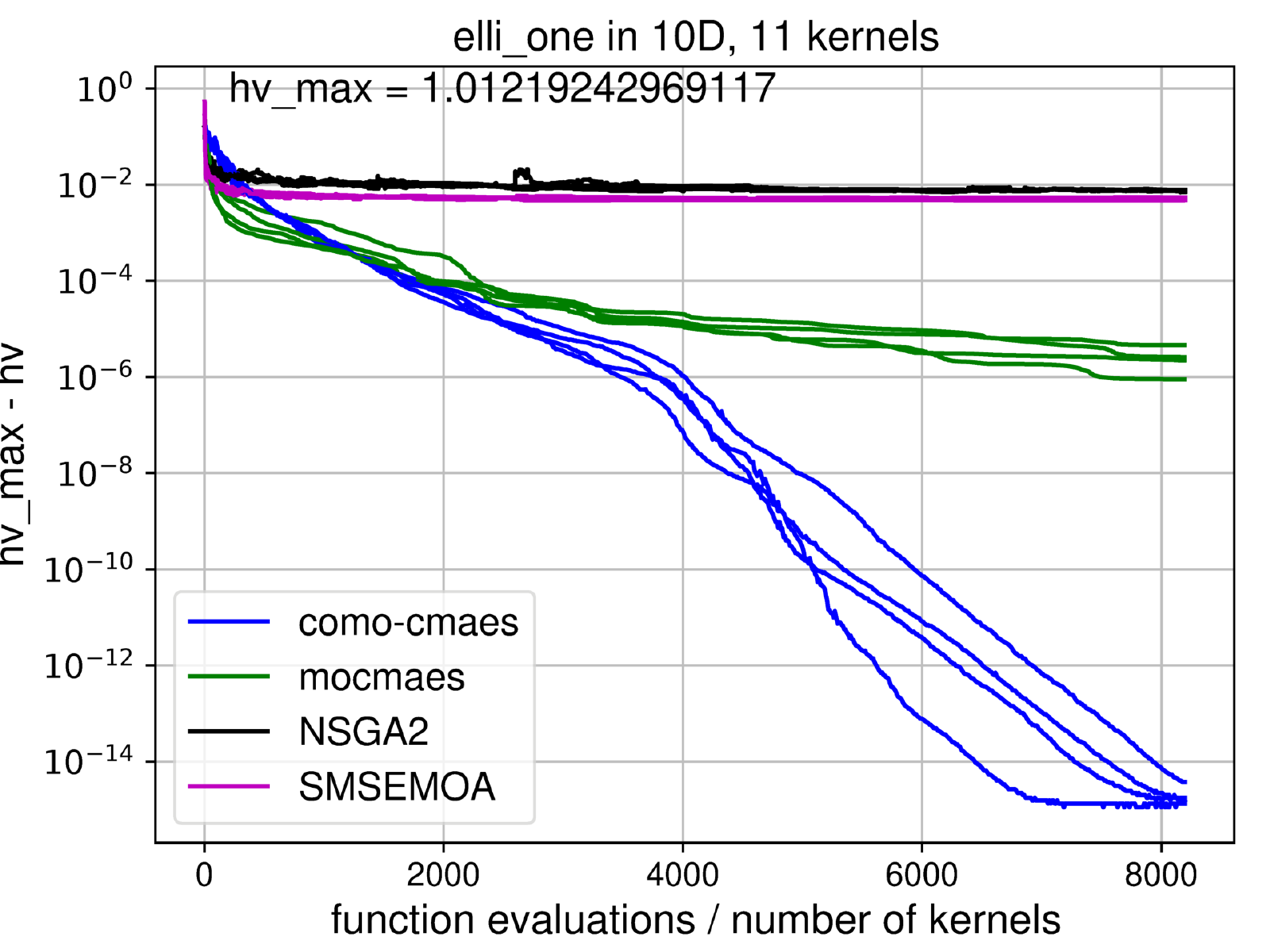}
	\includegraphics[width=0.24\textwidth]{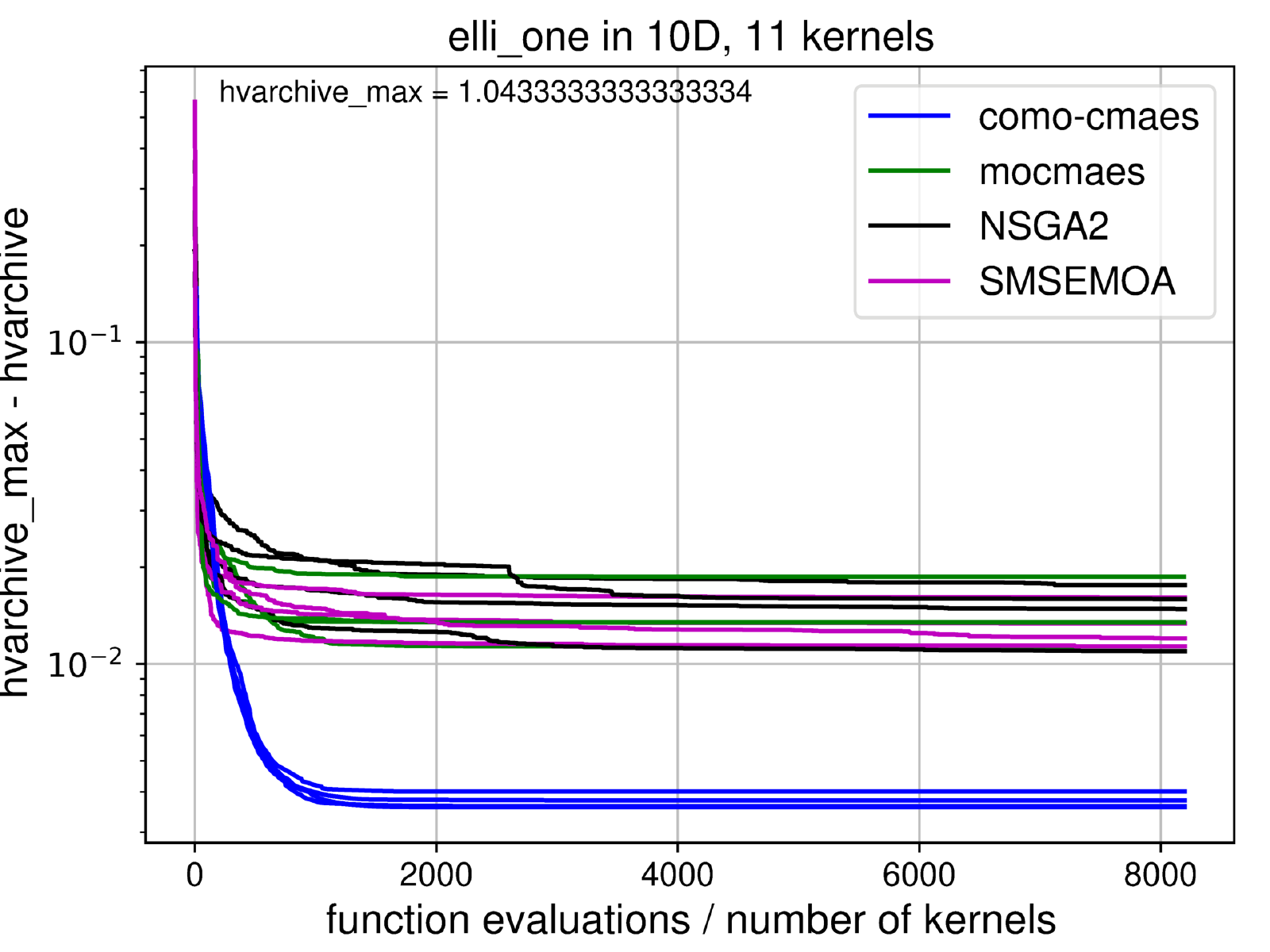}
	\includegraphics[width=0.24\textwidth]{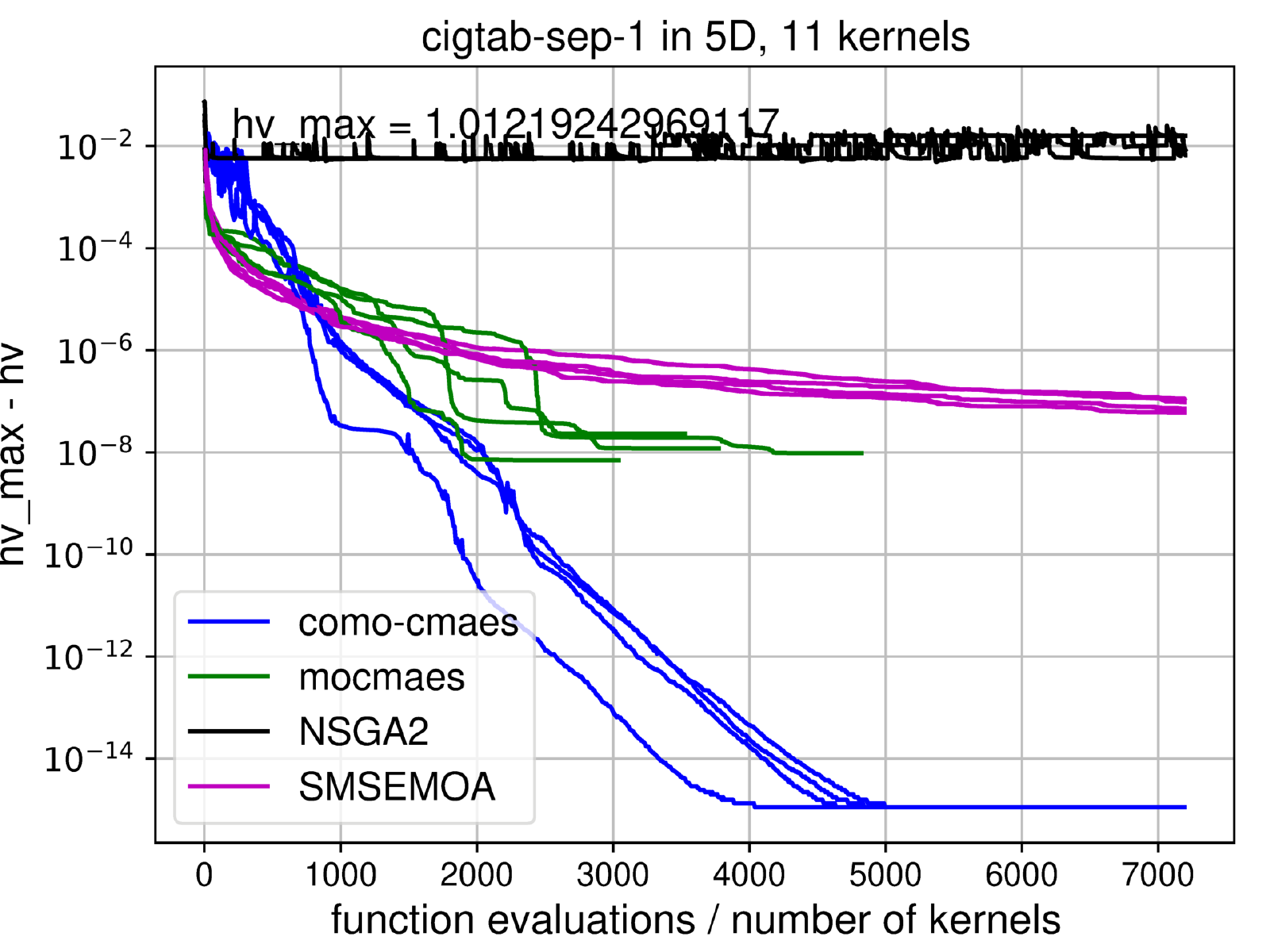}
	\includegraphics[width=0.24\textwidth]{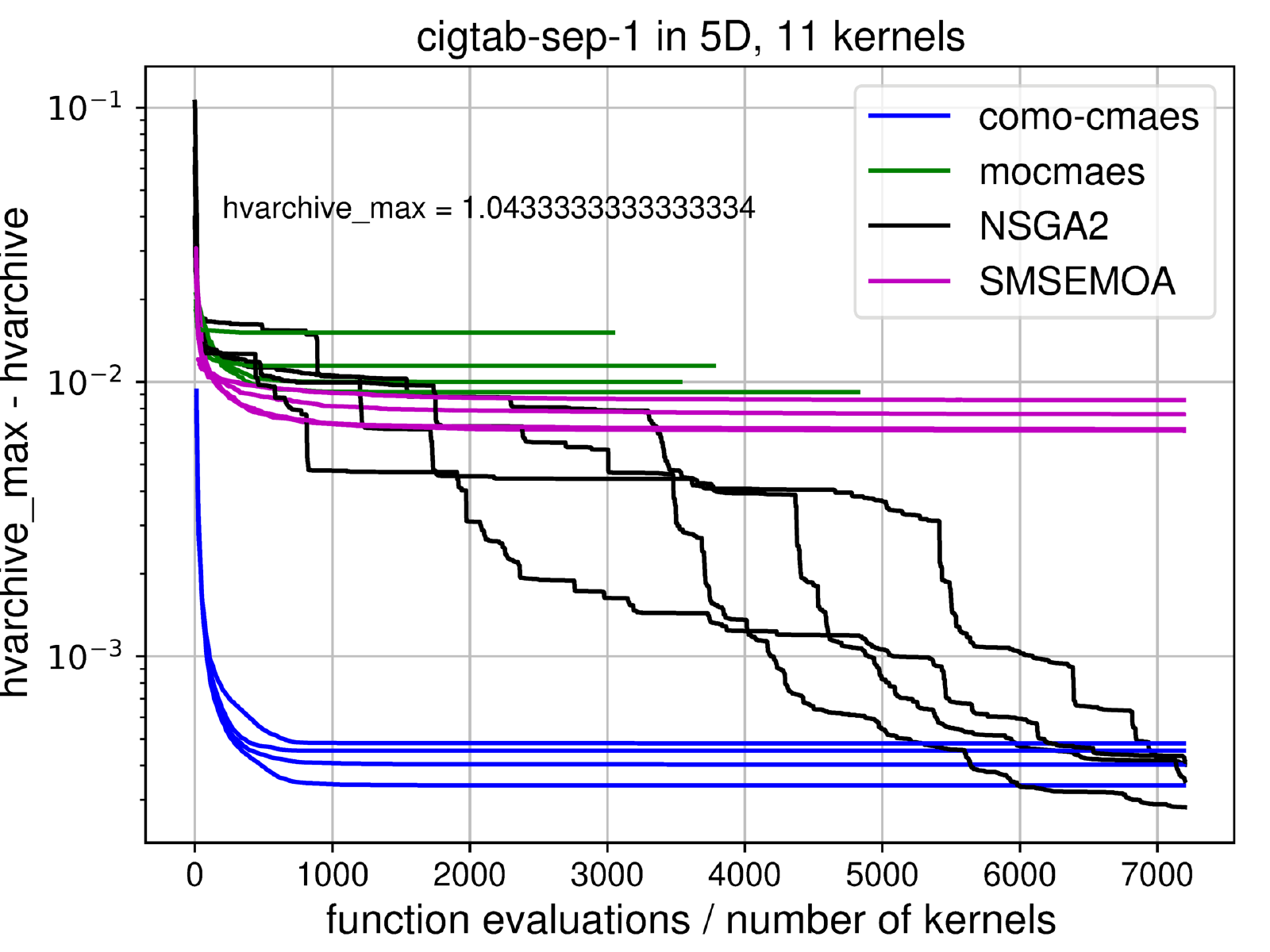}
	\includegraphics[width=0.24\textwidth]{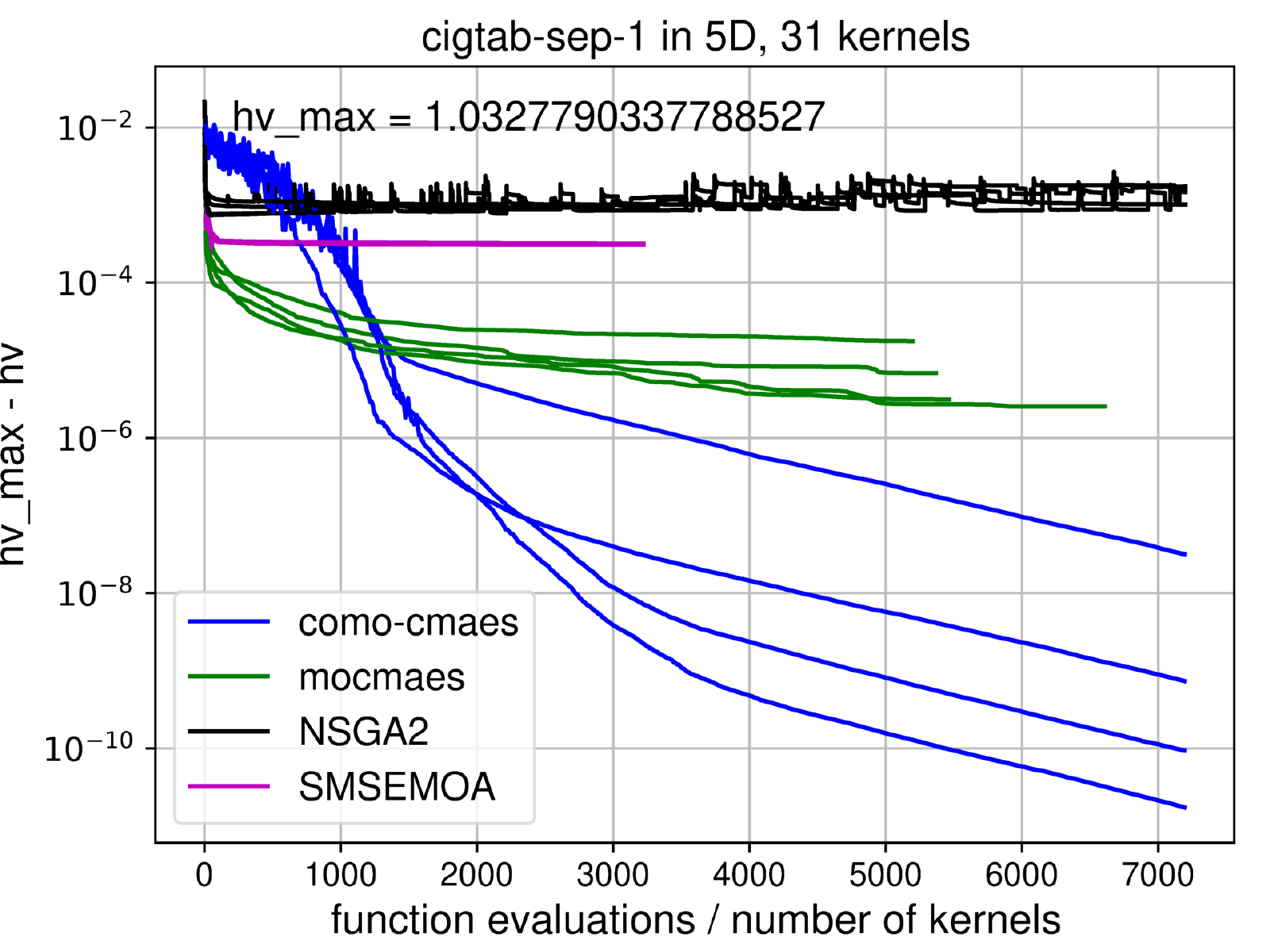}
	\includegraphics[width=0.24\textwidth]{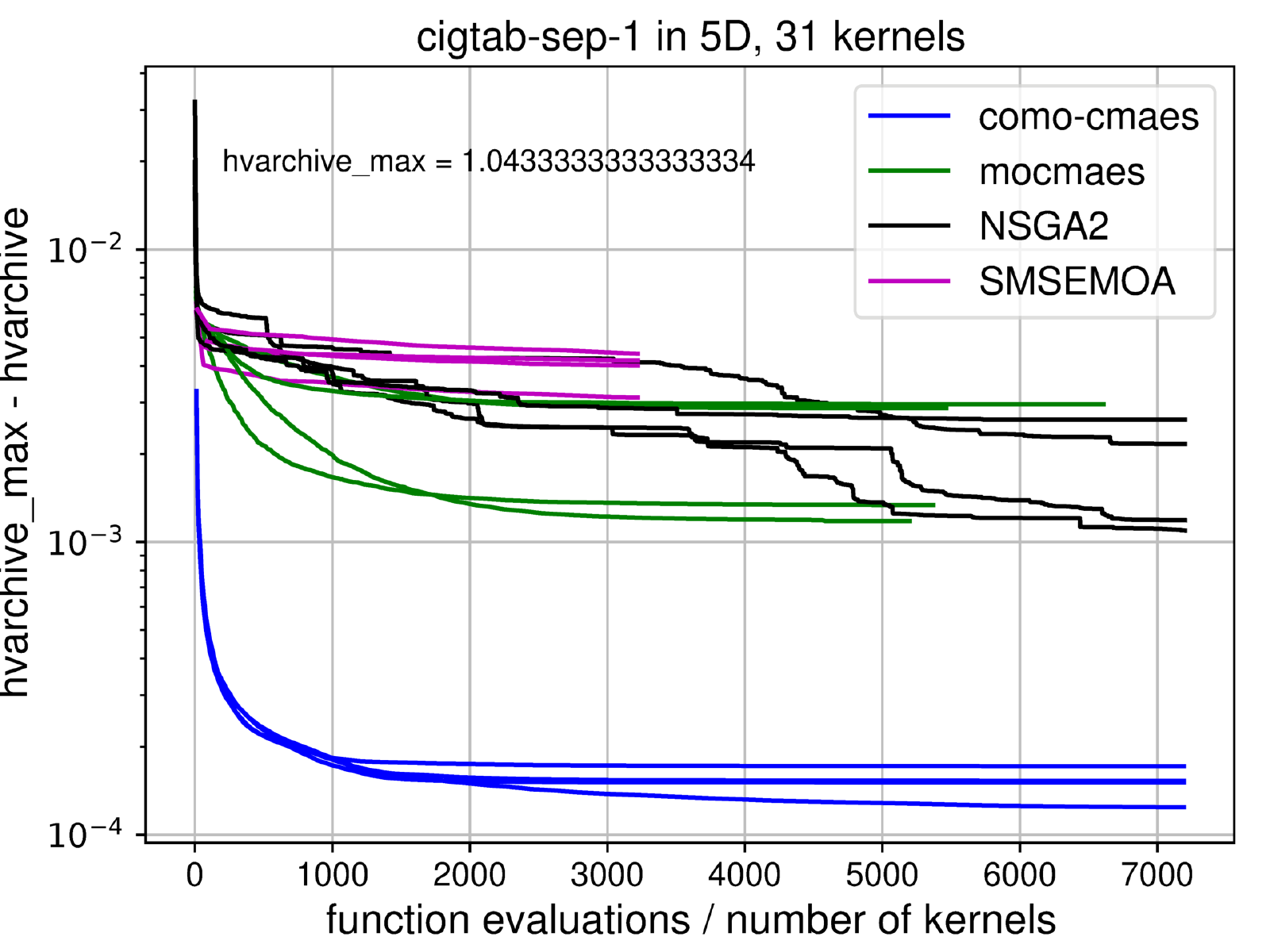}
	\includegraphics[width=0.24\textwidth]{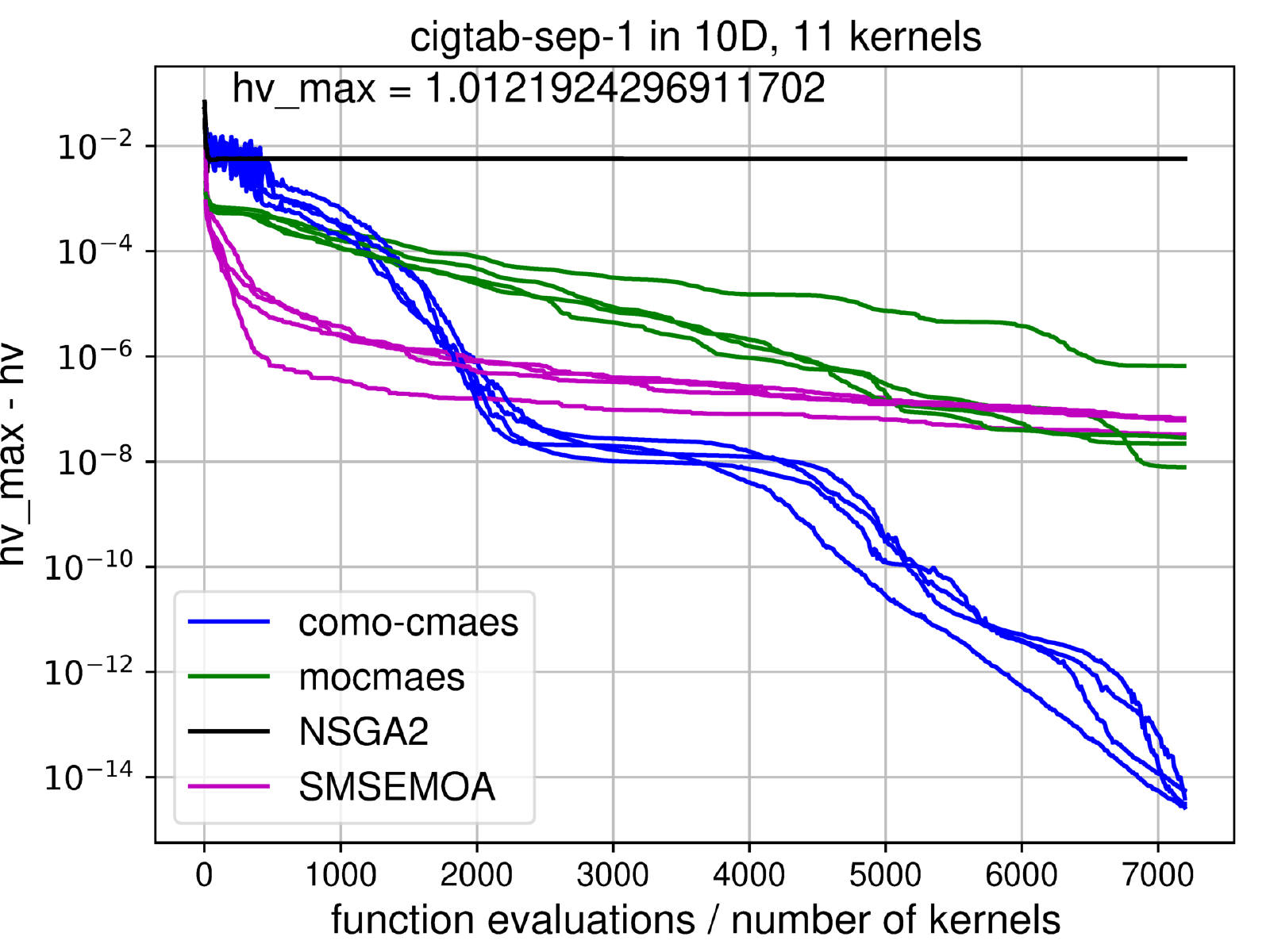}
	\includegraphics[width=0.24\textwidth]{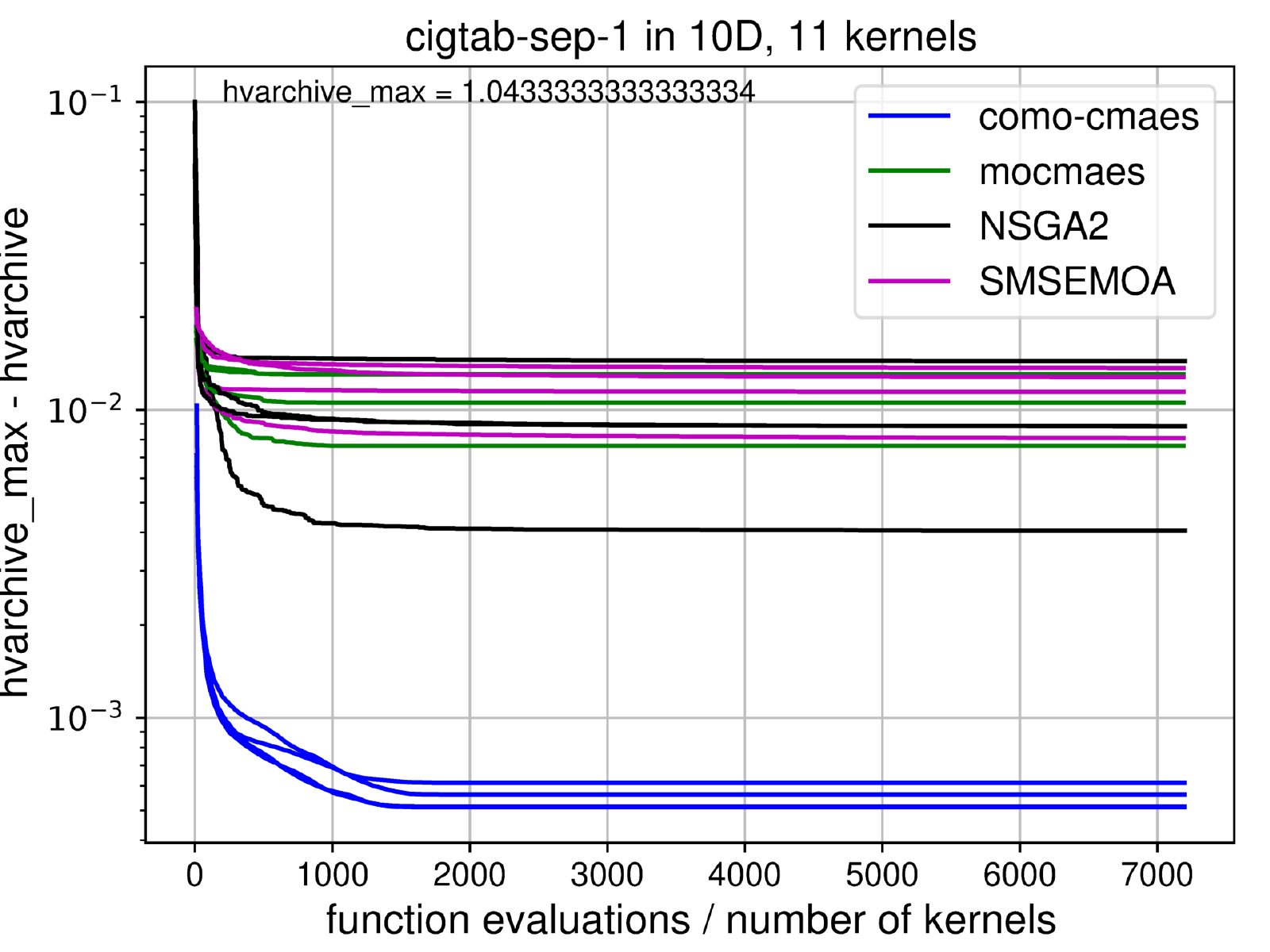}
	\includegraphics[width=0.24\textwidth]{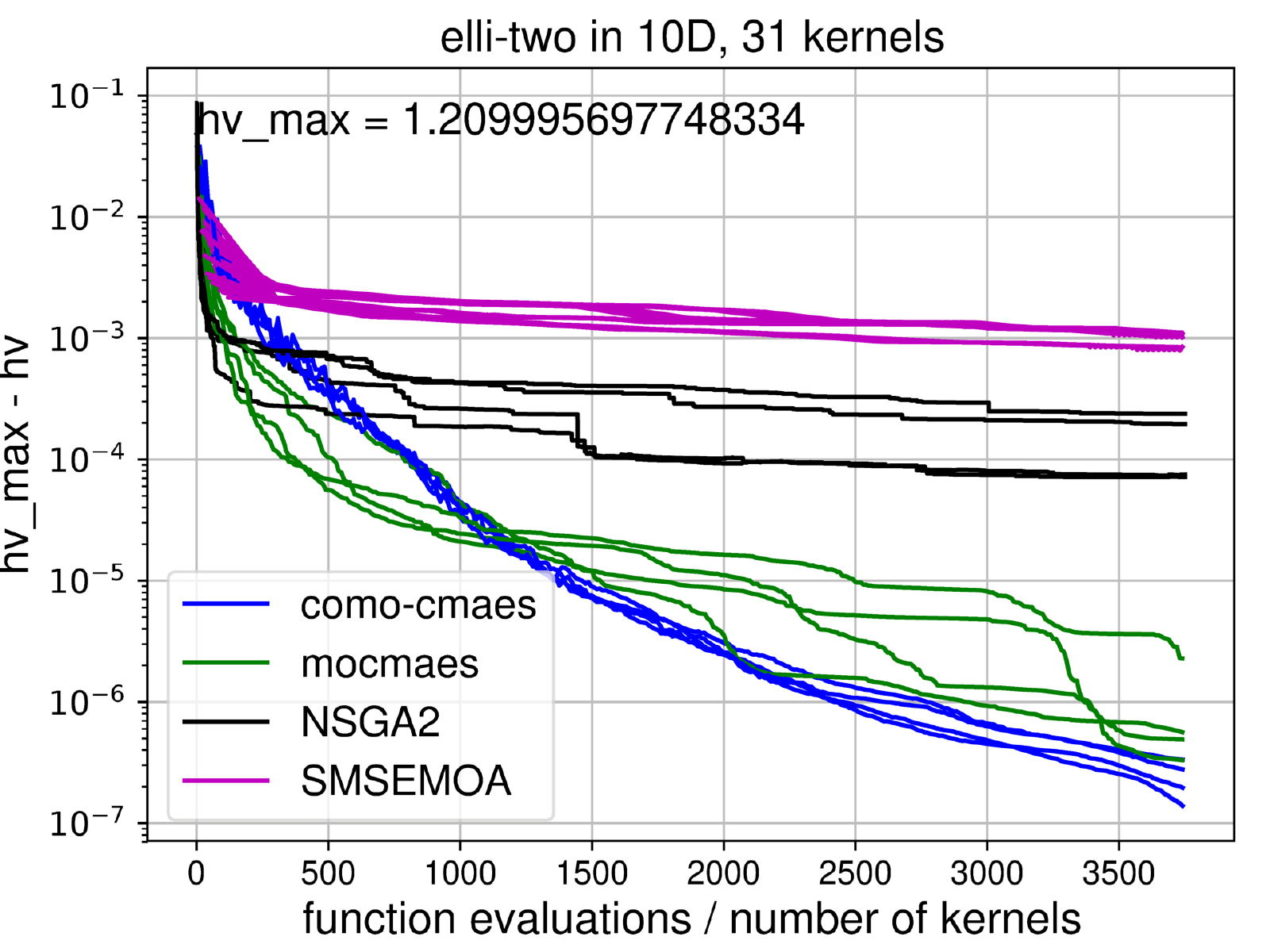}
	\includegraphics[width=0.24\textwidth]{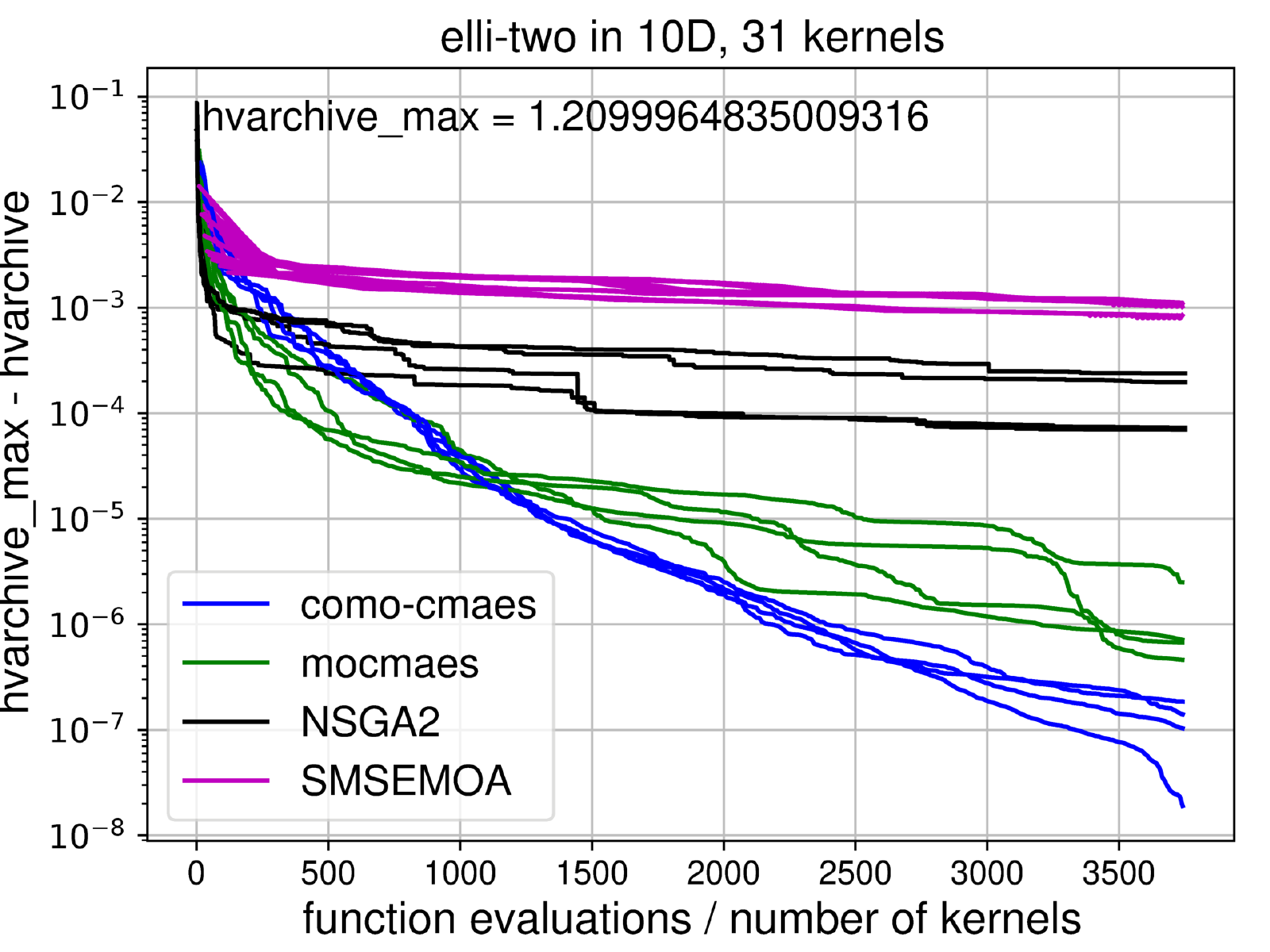}
	\caption{\label{fig:hypgap}  
	Convergence gaps (odd columns) and \archivegap{}s (even columns) for \textbf{bi-sphere}, \textbf{elli-one}, \textbf{cigtab-sep-$1$} and \textbf{elli-two}. Each algorithm is run $4$ times, in $5$D or $10$D, with $11$ or $31$ \kernels. The random matrices are drawn from the same seed in all the algorithms.}
\end{figure*}

\subsection{Comparing \COMOCMA  \,with MO-CMA-ES,  NSGA-II and SMS-EMOA}
We compare four multiobjective algorithms: COMO-CMA-ES, MO-CMA-ES \cite{ihr2007a}, NSGA-II \cite{dapm2002a} and SMS-EMOA \cite{bne2007a}, by testing them on classes of bi-objective convex-quadratic problems. We draw once and for all one rotation for \textbf{elli-one} in $10D$ and two different rotations for \textbf{elli-two} in $10D$. \del{ We run our problems with at most $17500\cdot \dimnumber$ function evaluations. }The Simulated Binary Crossover operator (SBX) and the polynomial mutation are used for NSGA-II (run with the \texttt{evoalgos} package~\cite{evoalgos}) and SMS-EMOA (run with the Matlab version by Fabian Kretzschmar and Tobias Wagner~\cite{wagner2010online}): we use a crossover probability of $0.7$ and a mutation probability of $0.1$, and the distribution indexes for crossover and mutation operators are both equal to $10$. We use the version of MO-CMA-ES from~\cite{vhi2010a}.
The number of \kernels for \COMOCMA corresponds to the population size of the other algorithms, that we set to either $11$ or $31$, and the dimensions considered are $5$ and $10$. The global initial step-size of \COMOCMA is set to $0.2$ with \nnew{initial} lower, upper bounds (line~8 of Algorithm~\ref{alg:como}) set to the all-zeros and all-ones vectors. The initial population for the three other algorithms is sampled uniformly at random in $[0,1]^n$.

We run each multiobjective optimization $4$ times and display the \convergencegap (of the population or the incumbent solutions of the \kernels) and the \archivegap.

In Figure~\ref{fig:hypgap}, the values of the \convergencegap reached by \COMOCMA and MO-CMA-ES are several orders of magnitude lower than for the two other algorithms.
\del{Only \COMOCMA and MO-CMA-ES achieve decent values in Figure~\ref{fig:hypgap} for the \convergencegap.}\del{ MO-CMA-ES reaches a precision between $10^{-12}$ and $10^{-14}$ for the \textbf{bi-sphere}, $5D$, $11$ kernels, with a better linear convergence rate than \COMOCMA that reaches a precision of $10^{-15}$.
}
On the 5-dimensional \textbf{bi-sphere}, \COMOCMA and MO-CMA-ES appear to show linear convergence, where the latter appears to be about 30\% faster than the former.
On the \textbf{cigtab-sep-$1$} function, \COMOCMA is initially slow, but catches up after about $1000$ evaluations per \kernel.
In all other cases, \COMOCMA\
shows superior performance for the \convergencegap.
On the 10-dimensional \textbf{cigtab-sep-$1$}, \COMOCMA\del{ with $11$ \kernels}
shows a plateau between 2000 and 4000 evaluations per \kernel.
This kind of plateau cannot be observed in the MO-CMA-ES and the observed final convergence speed is better for \COMOCMA than for MO-CMA-ES.
The observed plateau is typical for the behavior of non-elitist multi-recombinative CMA-ES on the tablet function, because CSA barely reduces an initially large step-size before the tablet-shape has been adapted,\niko{(confirmed with CSA vs TPA experiment)} which is related to the neutral subspace defect found in \cite{krause2017qualitative}. Elitism as in the MO-CMA-ES, on the other hand, also helps to decrease an initially too large step-size.\niko{ (confirmed with {\tt 'CMA\_elitist'='initial'} experiment)}

Although \COMOCMA was not designed to perform well on the \archivegap, it shows
consistently the best results over all experiments.
Only on the \textbf{cigtab-sep-$1$} in $5D$ with $11$ \kernels, NSGA-II reaches and slightly surpasses the \archivegap of \COMOCMA after $7000$ function evaluations per kernel.
This suggests, as expected from the known dependency between optimal step-size and population size \cite{hansen2015evolution}, that \COMOCMA adds valuable diversity while approaching the optimal $\kernumber$-distribution of the Pareto front at the same time.

\section{Conclusions}\label{sec:conclusions}
We have proposed (i) the \SUMO framework to define multiobjective optimizers from single-objective ones, (ii) a fitness for dominated solutions to be the distance to the empirical Pareto front (Uncrowded Hypervolume Improvement \UHVI) and (iii) the non-elitist "comma" CMA-ES to instantiate the framework (\COMOCMA). We observe that 
\COMOCMA converges linearly towards the $\kernumber$-optimal distribution of the hypervolume indicator on several bi-objective convex quadratic problems. The \COMOCMA  appears to be robust to independently rotating the Hessian matrices of convex-quadratic problems, even if such rotations transform the Pareto set from a line segment to a bent curve.
In our limited experiments, \COMOCMA performed generally better than MO-CMA-ES,  SMS-EMOA and the NSGA-II, w.r.t.\ \del{the }\convergencegap and \archivegap while \COMOCMA was solely designed to optimize the \convergencegap.
We conjecture that the advantage on the \archivegap\ is due to (i) the large stationary variance obtained with non-elitist evolution strategies and (ii) the fitness assignment of dominated solutions which favors the \nnew{vacant (uncrowded)} space between non-dominated solutions and hence serves as implicit crowding distance penalty measure.

%\noindent
%\\

%\noindent \paragraph{\textbf{% use this if space is an issue
\section*{% instead of this
% \textit{ Acknowledgements}}}
 Acknowledgements}
Part of this research has been conducted in the context of a research collaboration between Storengy and Inria. We particularly thank F. Huguet and A. Lange from Storengy for their strong support, practical ideas and expertise.

%\clearpage
%\newpage
%\balance
\bibliographystyle{ACM-Reference-Format}
\bibliography{gecco,allDimo}


\begin{thebibliography}{32}

%%% the Web address.

\ifx \showCODEN    \undefined \def \showCODEN     #1{\unskip}     \fi
\ifx \showDOI      \undefined \def \showDOI       #1{#1}\fi
\ifx \showISBNx    \undefined \def \showISBNx     #1{\unskip}     \fi
\ifx \showISBNxiii \undefined \def \showISBNxiii  #1{\unskip}     \fi
\ifx \showISSN     \undefined \def \showISSN      #1{\unskip}     \fi
\ifx \showLCCN     \undefined \def \showLCCN      #1{\unskip}     \fi
\ifx \shownote     \undefined \def \shownote      #1{#1}          \fi
\ifx \showarticletitle \undefined \def \showarticletitle #1{#1}   \fi
\ifx \showURL      \undefined \def \showURL       {\relax}        \fi
% invisible to TeX
\providecommand\bibfield[2]{#2}
\providecommand\bibinfo[2]{#2}
\providecommand\natexlab[1]{#1}
\providecommand\showeprint[2][]{arXiv:#2}

\bibitem[\protect\citeauthoryear{Auger, Bader, Brockhoff, and Zitzler}{Auger
  et~al\mbox{.}}{2009}]%
        {auger2009theory}
\bibfield{author}{\bibinfo{person}{Anne Auger}, \bibinfo{person}{Johannes
  Bader}, \bibinfo{person}{Dimo Brockhoff}, {and} \bibinfo{person}{Eckart
  Zitzler}.} \bibinfo{year}{2009}\natexlab{}.
\newblock \showarticletitle{Theory of the hypervolume indicator: optimal
  $\mu$-distributions and the choice of the reference point}. In
  \bibinfo{booktitle}{{\em Foundations of Genetic Algorithms (FOGA 2009)}}.
  \bibinfo{publisher}{ACM}, \bibinfo{address}{Orlando, Florida, USA},
  \bibinfo{pages}{87--102}.
\newblock


\bibitem[\protect\citeauthoryear{Auger, Bader, Brockhoff, and Zitzler}{Auger
  et~al\mbox{.}}{2012}]%
        {auger2012hypervolume}
\bibfield{author}{\bibinfo{person}{Anne Auger}, \bibinfo{person}{Johannes
  Bader}, \bibinfo{person}{Dimo Brockhoff}, {and} \bibinfo{person}{Eckart
  Zitzler}.} \bibinfo{year}{2012}\natexlab{}.
\newblock \showarticletitle{Hypervolume-based multiobjective optimization:
  Theoretical foundations and practical implications}.
\newblock \bibinfo{journal}{{\em Theoretical Computer Science\/}}
  \bibinfo{volume}{425} (\bibinfo{year}{2012}), \bibinfo{pages}{75--103}.
\newblock


\bibitem[\protect\citeauthoryear{Berghammer, Friedrich, and Neumann}{Berghammer
  et~al\mbox{.}}{2010}]%
        {bfn2010a}
\bibfield{author}{\bibinfo{person}{R. Berghammer}, \bibinfo{person}{T.
  Friedrich}, {and} \bibinfo{person}{F. Neumann}.}
  \bibinfo{year}{2010}\natexlab{}.
\newblock \showarticletitle{{Set-based Multi-objective Optimization,
  Indicators, and Deteriorative Cycles}}. In \bibinfo{booktitle}{{\em Genetic
  and Evolutionary Computation Conference (GECCO 2010)}}.
  \bibinfo{publisher}{ACM}, \bibinfo{address}{Portland, Oregon},
  \bibinfo{pages}{495--502}.
\newblock
\showDOI{%
\url{https://doi.org/10.1145/1830483.1830574}}


\bibitem[\protect\citeauthoryear{Beume, Naujoks, and Emmerich}{Beume
  et~al\mbox{.}}{2007}]%
        {bne2007a}
\bibfield{author}{\bibinfo{person}{N. Beume}, \bibinfo{person}{B. Naujoks},
  {and} \bibinfo{person}{M. Emmerich}.} \bibinfo{year}{2007}\natexlab{}.
\newblock \showarticletitle{{SMS-EMOA: Multiobjective Selection Based on
  Dominated Hypervolume}}.
\newblock \bibinfo{journal}{{\em European Journal of Operational Research\/}}
  \bibinfo{volume}{181}, \bibinfo{number}{3} (\bibinfo{year}{2007}),
  \bibinfo{pages}{1653--1669}.
\newblock


\bibitem[\protect\citeauthoryear{Bringmann and Friedrich}{Bringmann and
  Friedrich}{2011}]%
        {bringmann2011convergence}
\bibfield{author}{\bibinfo{person}{Karl Bringmann} {and}
  \bibinfo{person}{Tobias Friedrich}.} \bibinfo{year}{2011}\natexlab{}.
\newblock \showarticletitle{Convergence of hypervolume-based archiving
  algorithms I: Effectiveness}. In \bibinfo{booktitle}{{\em Proceedings of the
  13th annual conference on Genetic and evolutionary computation}}.
  \bibinfo{publisher}{ACM}, \bibinfo{address}{Dublin, Ireland},
  \bibinfo{pages}{745--752}.
\newblock


\bibitem[\protect\citeauthoryear{Bubeck and Cesa-Bianchi}{Bubeck and
  Cesa-Bianchi}{2012}]%
        {bc2012a}
\bibfield{author}{\bibinfo{person}{S{\'e}bastien Bubeck} {and}
  \bibinfo{person}{Nicolo Cesa-Bianchi}.} \bibinfo{year}{2012}\natexlab{}.
\newblock \showarticletitle{Regret analysis of stochastic and nonstochastic
  multi-armed bandit problems}.
\newblock \bibinfo{journal}{{\em Foundations and Trends{\textregistered} in
  Machine Learning\/}} \bibinfo{volume}{5}, \bibinfo{number}{1}
  (\bibinfo{year}{2012}), \bibinfo{pages}{1--122}.
\newblock


\bibitem[\protect\citeauthoryear{Collette, Hansen, Pujol, Salazar~Aponte, and
  Le~Riche}{Collette et~al\mbox{.}}{2010}]%
        {collette2010object}
\bibfield{author}{\bibinfo{person}{Yann Collette}, \bibinfo{person}{Nikolaus
  Hansen}, \bibinfo{person}{Gilles Pujol}, \bibinfo{person}{Daniel
  Salazar~Aponte}, {and} \bibinfo{person}{Rodolphe Le~Riche}.}
  \bibinfo{year}{2010}\natexlab{}.
\newblock \showarticletitle{On Object-Oriented Programming of Optimizers -
  Examples in Scilab}.
\newblock In \bibinfo{booktitle}{{\em Multidisciplinary Design Optimization in
  Computational Mechanics}}, \bibfield{editor}{\bibinfo{person}{Rajan~Filomeno
  Coelho} {and} \bibinfo{person}{Piotr Breitkopf}} (Eds.).
  \bibinfo{publisher}{Wiley}, \bibinfo{address}{New Jersey},
  \bibinfo{pages}{499--538}.
\newblock
\showURL{%
\url{https://hal.inria.fr/inria-00476172}}


\bibitem[\protect\citeauthoryear{Deb, Pratap, Agarwal, and Meyarivan}{Deb
  et~al\mbox{.}}{2002}]%
        {dapm2002a}
\bibfield{author}{\bibinfo{person}{K. Deb}, \bibinfo{person}{A. Pratap},
  \bibinfo{person}{S. Agarwal}, {and} \bibinfo{person}{T. Meyarivan}.}
  \bibinfo{year}{2002}\natexlab{}.
\newblock \showarticletitle{{A Fast and Elitist Multiobjective Genetic
  Algorithm: NSGA-II}}.
\newblock \bibinfo{journal}{{\em IEEE Transactions on Evolutionary
  Computation\/}} \bibinfo{volume}{6}, \bibinfo{number}{2}
  (\bibinfo{year}{2002}), \bibinfo{pages}{182--197}.
\newblock


\bibitem[\protect\citeauthoryear{Emmerich, Beume, and Naujoks}{Emmerich
  et~al\mbox{.}}{2005}]%
        {emmerich2005emo}
\bibfield{author}{\bibinfo{person}{Michael Emmerich}, \bibinfo{person}{Nicola
  Beume}, {and} \bibinfo{person}{Boris Naujoks}.}
  \bibinfo{year}{2005}\natexlab{}.
\newblock \showarticletitle{An EMO algorithm using the hypervolume measure as
  selection criterion}. In \bibinfo{booktitle}{{\em International Conference on
  Evolutionary Multi-Criterion Optimization}}. \bibinfo{publisher}{Springer},
  \bibinfo{address}{Guanajuato, Mexico}, \bibinfo{pages}{62--76}.
\newblock


\bibitem[\protect\citeauthoryear{Emmerich and Klinkenberg}{Emmerich and
  Klinkenberg}{2008}]%
        {ek2008a}
\bibfield{author}{\bibinfo{person}{Michael Emmerich} {and}
  \bibinfo{person}{Jan-willem Klinkenberg}.} \bibinfo{year}{2008}\natexlab{}.
\newblock \bibinfo{booktitle}{{\em The computation of the expected improvement
  in dominated hypervolume of Pareto front approximations}}.
\newblock \bibinfo{type}{Technical Report} 4-2008. \bibinfo{institution}{Leiden
  Institute of Advanced Computer Science, LIACS}.
\newblock


\bibitem[\protect\citeauthoryear{Goldberg}{Goldberg}{1989}]%
        {gold1989a}
\bibfield{author}{\bibinfo{person}{D.~E. Goldberg}.}
  \bibinfo{year}{1989}\natexlab{}.
\newblock \bibinfo{booktitle}{{\em {Genetic Algorithms in Search, Optimization,
  and Machine Learning}}}.
\newblock \bibinfo{publisher}{Addison-Wesley}, \bibinfo{address}{Reading,
  Massachusetts}.
\newblock


\bibitem[\protect\citeauthoryear{Hansen and Jaszkiewicz}{Hansen and
  Jaszkiewicz}{1998}]%
        {hj1998a}
\bibfield{author}{\bibinfo{person}{M.~P. Hansen} {and} \bibinfo{person}{A.
  Jaszkiewicz}.} \bibinfo{year}{1998}\natexlab{}.
\newblock \bibinfo{booktitle}{{\em {Evaluating The Quality of Approximations of
  the Non-Dominated Set}}}.
\newblock \bibinfo{type}{{T}echnical {R}eport}. \bibinfo{institution}{Institute
  of Mathematical Modeling, Technical University of Denmark}.
\newblock
\newblock
\shownote{IMM Technical Report IMM-REP-1998-7.}


\bibitem[\protect\citeauthoryear{Hansen, Akimoto, and Baudis}{Hansen
  et~al\mbox{.}}{2019}]%
        {hansen2019pycma}
\bibfield{author}{\bibinfo{person}{Nikolaus Hansen}, \bibinfo{person}{Youhei
  Akimoto}, {and} \bibinfo{person}{Petr Baudis}.}
  \bibinfo{year}{2019}\natexlab{}.
\newblock \bibinfo{title}{{CMA-ES/pycma} on {G}ithub}.
\newblock \bibinfo{howpublished}{Zenodo, DOI:10.5281/zenodo.2559634}.
  (\bibinfo{date}{Feb.} \bibinfo{year}{2019}).
\newblock
\showDOI{%
\url{https://doi.org/10.5281/zenodo.2559634}}


\bibitem[\protect\citeauthoryear{Hansen, Arnold, and Auger}{Hansen
  et~al\mbox{.}}{2015}]%
        {hansen2015evolution}
\bibfield{author}{\bibinfo{person}{Nikolaus Hansen}, \bibinfo{person}{Dirk~V
  Arnold}, {and} \bibinfo{person}{Anne Auger}.}
  \bibinfo{year}{2015}\natexlab{}.
\newblock \showarticletitle{Evolution strategies}.
\newblock In \bibinfo{booktitle}{{\em Springer handbook of computational
  intelligence}}. \bibinfo{publisher}{Springer}, \bibinfo{address}{Berlin},
  \bibinfo{pages}{871--898}.
\newblock


\bibitem[\protect\citeauthoryear{Hansen and Ostermeier}{Hansen and
  Ostermeier}{2001}]%
        {ho2001a}
\bibfield{author}{\bibinfo{person}{N. Hansen} {and} \bibinfo{person}{A.
  Ostermeier}.} \bibinfo{year}{2001}\natexlab{}.
\newblock \showarticletitle{{Completely Derandomized Self-Adaptation in
  Evolution Strategies}}.
\newblock \bibinfo{journal}{{\em Evolutionary Computation\/}}
  \bibinfo{volume}{9}, \bibinfo{number}{2} (\bibinfo{year}{2001}),
  \bibinfo{pages}{159--195}.
\newblock


\bibitem[\protect\citeauthoryear{Hernandez, Schutze, Wang, Deutz, and
  Emmerich}{Hernandez et~al\mbox{.}}{2018}]%
        {hernandez2018set}
\bibfield{author}{\bibinfo{person}{VAS Hernandez}, \bibinfo{person}{O Schutze},
  \bibinfo{person}{H Wang}, \bibinfo{person}{A Deutz}, {and} \bibinfo{person}{M
  Emmerich}.} \bibinfo{year}{2018}\natexlab{}.
\newblock \showarticletitle{The Set-Based Hypervolume Newton Method for
  Bi-Objective Optimization}.
\newblock \bibinfo{journal}{{\em IEEE transactions on cybernetics\/}}
  \bibinfo{volume}{in print} (\bibinfo{year}{2018}).
\newblock
\newblock
\shownote{(in print).}


\bibitem[\protect\citeauthoryear{Igel, Hansen, and Roth}{Igel
  et~al\mbox{.}}{2007}]%
        {ihr2007a}
\bibfield{author}{\bibinfo{person}{C. Igel}, \bibinfo{person}{N. Hansen}, {and}
  \bibinfo{person}{S. Roth}.} \bibinfo{year}{2007}\natexlab{}.
\newblock \showarticletitle{Covariance matrix adaptation for multi-objective
  optimization}.
\newblock \bibinfo{journal}{{\em Evolutionary Computation\/}}
  \bibinfo{volume}{15}, \bibinfo{number}{1} (\bibinfo{year}{2007}),
  \bibinfo{pages}{1--28}.
\newblock


\bibitem[\protect\citeauthoryear{Keane}{Keane}{2006}]%
        {kean2006a}
\bibfield{author}{\bibinfo{person}{{Andy J.} Keane}.}
  \bibinfo{year}{2006}\natexlab{}.
\newblock \showarticletitle{Statistical improvement criteria for use in
  multiobjective design optimization}.
\newblock \bibinfo{journal}{{\em AIAA journal\/}} \bibinfo{volume}{44},
  \bibinfo{number}{4} (\bibinfo{year}{2006}), \bibinfo{pages}{879--891}.
\newblock


\bibitem[\protect\citeauthoryear{Knowles, Thiele, and Zitzler}{Knowles
  et~al\mbox{.}}{2006}]%
        {ktz2006a}
\bibfield{author}{\bibinfo{person}{J. Knowles}, \bibinfo{person}{L. Thiele},
  {and} \bibinfo{person}{E. Zitzler}.} \bibinfo{year}{2006}\natexlab{}.
\newblock \bibinfo{booktitle}{{\em {A Tutorial on the Performance Assessment of
  Stochastic Multiobjective Optimizers}}}.
\newblock \bibinfo{type}{TIK Report} 214. \bibinfo{institution}{Computer
  Engineering and Networks Laboratory (TIK), ETH Zurich}.
\newblock


\bibitem[\protect\citeauthoryear{Krause, Glasmachers, and Igel}{Krause
  et~al\mbox{.}}{2017}]%
        {krause2017qualitative}
\bibfield{author}{\bibinfo{person}{Oswin Krause}, \bibinfo{person}{Tobias
  Glasmachers}, {and} \bibinfo{person}{Christian Igel}.}
  \bibinfo{year}{2017}\natexlab{}.
\newblock \showarticletitle{Qualitative and quantitative assessment of step
  size adaptation rules}. In \bibinfo{booktitle}{{\em Proceedings of the 14th
  ACM/SIGEVO Conference on Foundations of Genetic Algorithms}}.
  \bibinfo{publisher}{ACM}, \bibinfo{address}{Copenhagen, Denmark},
  \bibinfo{pages}{139--148}.
\newblock


\bibitem[\protect\citeauthoryear{Miettinen}{Miettinen}{1999}]%
        {miet1999a}
\bibfield{author}{\bibinfo{person}{K. Miettinen}.}
  \bibinfo{year}{1999}\natexlab{}.
\newblock \bibinfo{booktitle}{{\em {Nonlinear Multiobjective Optimization}}}.
\newblock \bibinfo{publisher}{Kluwer}, \bibinfo{address}{Boston, MA, USA}.
\newblock


\bibitem[\protect\citeauthoryear{Ponweiser, Wagner, Biermann, and
  Vincze}{Ponweiser et~al\mbox{.}}{2008}]%
        {pwbv2008a}
\bibfield{author}{\bibinfo{person}{Wolfgang Ponweiser}, \bibinfo{person}{Tobias
  Wagner}, \bibinfo{person}{Dirk Biermann}, {and} \bibinfo{person}{Markus
  Vincze}.} \bibinfo{year}{2008}\natexlab{}.
\newblock \showarticletitle{Multiobjective Optimization on a Limited Budget of
  Evaluations Using Model-Assisted ${\mathcal{S}}$-Metric Selection}. In
  \bibinfo{booktitle}{{\em Parallel Problem Solving from Nature (PPSN 2008)}}.
  \bibinfo{publisher}{Springer}, \bibinfo{address}{Dortmund, Germany},
  \bibinfo{pages}{784--794}.
\newblock


\bibitem[\protect\citeauthoryear{Toure, Auger, Brockhoff, and Hansen}{Toure
  et~al\mbox{.}}{2019}]%
        {toure2019bi}
\bibfield{author}{\bibinfo{person}{Cheikh Toure}, \bibinfo{person}{Anne Auger},
  \bibinfo{person}{Dimo Brockhoff}, {and} \bibinfo{person}{Nikolaus Hansen}.}
  \bibinfo{year}{2019}\natexlab{}.
\newblock \showarticletitle{On Bi-Objective convex-quadratic problems}. In
  \bibinfo{booktitle}{{\em International Conference on Evolutionary
  Multi-Criterion Optimization}}. \bibinfo{publisher}{Springer},
  \bibinfo{address}{Lansing, Michigan, USA}, \bibinfo{pages}{3--14}.
\newblock


\bibitem[\protect\citeauthoryear{Vo\ss, Hansen, and Igel}{Vo\ss
  et~al\mbox{.}}{2010}]%
        {vhi2010a}
\bibfield{author}{\bibinfo{person}{T. Vo\ss}, \bibinfo{person}{N. Hansen},
  {and} \bibinfo{person}{C. Igel}.} \bibinfo{year}{2010}\natexlab{}.
\newblock \showarticletitle{{Improved Step Size Adaptation for the MO-CMA-ES}}.
  In \bibinfo{booktitle}{{\em Genetic and Evolutionary Computation Conference
  (GECCO 2010)}}, \bibfield{editor}{\bibinfo{person}{J.~Branke}
  {et~al\mbox{.}}} (Eds.). \bibinfo{publisher}{ACM},
  \bibinfo{address}{Portland, OR, USA}, \bibinfo{pages}{487--494}.
\newblock


\bibitem[\protect\citeauthoryear{Wagner, Emmerich, Deutz, and Ponweiser}{Wagner
  et~al\mbox{.}}{2010}]%
        {wedp2010a}
\bibfield{author}{\bibinfo{person}{Tobias Wagner}, \bibinfo{person}{Michael
  Emmerich}, \bibinfo{person}{Andr{\'e} Deutz}, {and} \bibinfo{person}{Wolfgang
  Ponweiser}.} \bibinfo{year}{2010}\natexlab{}.
\newblock \showarticletitle{On expected-improvement criteria for model-based
  multi-objective optimization}. In \bibinfo{booktitle}{{\em International
  Conference on Parallel Problem Solving from Nature}}.
  \bibinfo{publisher}{Springer}, \bibinfo{address}{Krakow, Poland},
  \bibinfo{pages}{718--727}.
\newblock


\bibitem[\protect\citeauthoryear{Wagner and Trautmann}{Wagner and
  Trautmann}{2010}]%
        {wagner2010online}
\bibfield{author}{\bibinfo{person}{Tobias Wagner} {and} \bibinfo{person}{Heike
  Trautmann}.} \bibinfo{year}{2010}\natexlab{}.
\newblock \showarticletitle{Online convergence detection for evolutionary
  multi-objective algorithms revisited}. In \bibinfo{booktitle}{{\em IEEE
  Congress on Evolutionary Computation}}. \bibinfo{publisher}{IEEE},
  \bibinfo{address}{Barcelona, Spain}, \bibinfo{pages}{1--8}.
\newblock


\bibitem[\protect\citeauthoryear{Wessing}{Wessing}{2017}]%
        {evoalgos}
\bibfield{author}{\bibinfo{person}{Simon Wessing}.}
  \bibinfo{year}{2017}\natexlab{}.
\newblock \bibinfo{title}{evoalgos: Modular evolutionary algorithms. Python
  package version 1}.
\newblock   (\bibinfo{year}{2017}).
\newblock
\showURL{%
\url{https://pypi.python.org/pypi/evoalgos}}
\newblock
\shownote{[Online; accessed 31-January-2019].}


\bibitem[\protect\citeauthoryear{Yang, Emmerich, Deutz, and B{\"a}ck}{Yang
  et~al\mbox{.}}{2019}]%
        {yang2019multi}
\bibfield{author}{\bibinfo{person}{Kaifeng Yang}, \bibinfo{person}{Michael
  Emmerich}, \bibinfo{person}{Andr{\'e} Deutz}, {and} \bibinfo{person}{Thomas
  B{\"a}ck}.} \bibinfo{year}{2019}\natexlab{}.
\newblock \showarticletitle{Multi-Objective Bayesian Global Optimization using
  expected hypervolume improvement gradient}.
\newblock \bibinfo{journal}{{\em Swarm and evolutionary computation\/}}
  \bibinfo{volume}{44} (\bibinfo{year}{2019}), \bibinfo{pages}{945--956}.
\newblock


\bibitem[\protect\citeauthoryear{Zhang and Li}{Zhang and Li}{2007}]%
        {zl2007a}
\bibfield{author}{\bibinfo{person}{Q. Zhang} {and} \bibinfo{person}{H. Li}.}
  \bibinfo{year}{2007}\natexlab{}.
\newblock \showarticletitle{{MOEA/D: A Multiobjective Evolutionary Algorithm
  Based on Decomposition}}.
\newblock \bibinfo{journal}{{\em IEEE Transactions on Evolutionary
  Computation\/}} \bibinfo{volume}{11}, \bibinfo{number}{6}
  (\bibinfo{year}{2007}), \bibinfo{pages}{712--731}.
\newblock
\showDOI{%
\url{https://doi.org/10.1109/TEVC.2007.892759}}


\bibitem[\protect\citeauthoryear{Zitzler and K{\"u}nzli}{Zitzler and
  K{\"u}nzli}{2004}]%
        {zitzler2004indicator}
\bibfield{author}{\bibinfo{person}{Eckart Zitzler} {and} \bibinfo{person}{Simon
  K{\"u}nzli}.} \bibinfo{year}{2004}\natexlab{}.
\newblock \showarticletitle{Indicator-based selection in multiobjective
  search}. In \bibinfo{booktitle}{{\em International Conference on Parallel
  Problem Solving from Nature}}. \bibinfo{publisher}{Springer},
  \bibinfo{address}{Birmingham, UK}, \bibinfo{pages}{832--842}.
\newblock


\bibitem[\protect\citeauthoryear{Zitzler and Thiele}{Zitzler and
  Thiele}{1998a}]%
        {zt1998b}
\bibfield{author}{\bibinfo{person}{E. Zitzler} {and} \bibinfo{person}{L.
  Thiele}.} \bibinfo{year}{1998}\natexlab{a}.
\newblock \showarticletitle{{Multiobjective Optimization Using Evolutionary
  Algorithms - A Comparative Case Study}}. In \bibinfo{booktitle}{{\em
  Conference on Parallel Problem Solving from Nature (PPSN~V)}} {\em
  (\bibinfo{series}{LNCS})}, Vol.~\bibinfo{volume}{1498}.
  \bibinfo{publisher}{Springer}, \bibinfo{address}{Amsterdam, The Netherlands},
  \bibinfo{pages}{292--301}.
\newblock


\bibitem[\protect\citeauthoryear{Zitzler and Thiele}{Zitzler and
  Thiele}{1998b}]%
        {zitzler1998multiobjective}
\bibfield{author}{\bibinfo{person}{Eckart Zitzler} {and}
  \bibinfo{person}{Lothar Thiele}.} \bibinfo{year}{1998}\natexlab{b}.
\newblock \showarticletitle{Multiobjective optimization using evolutionary
  algorithms - A comparative case study}. In \bibinfo{booktitle}{{\em
  International conference on parallel problem solving from nature}}.
  \bibinfo{publisher}{Springer}, \bibinfo{address}{Amsterdam, The Netherlands},
  \bibinfo{pages}{292--301}.
\newblock


\end{thebibliography}

%\end{document}
%
%
%
%\pagebreak
%
%\appendix
%\section{WHAT FOLLOWS IS OLD STUFF (which we still might want to integrate}
%
%
%\todo{Dimo}
%\begin{itemize}
%\item [Dimo] Multiobjective optimization is important, look for the Pareto-set instead of a single point, typically it is a curve or union of curve [check that curve is really the term]. For instance, double sphere, segment between both optima.
%\item [Dimo] In this context, convergence to the Pareto-front is not immediate to define.
%\item At the same time, key of success of single-objective algorithms like CMA-ES is that they exhibit linear convergence on wide classes of functions (+ learning second order information).
%\item Natural to ask that this linear convergence is transferred in the MO context.
%\item Contributions of the paper: (1) propose a general framework for building a MO algorithm from a single-objective (SO) one ; composed of $\objnumber$ instances of the SO algorithm where each converge linearly to a point of the Pareto set.
%\item (2) Instantiation of the framework using the CMA-ES algorithm, naturally obtained a new MO algorithm, the COMO-CMA-ES.
%\item (3) Empirical validation ; illustration of the linear convergence
%\end{itemize}

\end{document}